\documentclass{article}%

\usepackage{amsmath,amsfonts,amssymb}
\usepackage{graphicx}%

\providecommand{\U}[1]{\protect\rule{.1in}{.1in}}

\newtheorem{theorem}{Theorem}
\newtheorem{algorithm}{Algorithm}
\newtheorem{corollary}[theorem]{Corollary}
\newtheorem{definition}{Definition}
\newtheorem{example}{Example}
\newtheorem{lemma}[theorem]{Lemma}
\newtheorem{remark}{Remark}
\numberwithin{equation}{section}

\newenvironment{proof}[1][Proof]{\noindent\textbf{#1.} }{\ \rule{0.5em}{0.5em}}

\begin{document}

\title{Explicit Local Time-Stepping for the Inhomogeneous
Wave Equation with Optimal Convergence}
\author{Marcus J. Grote\thanks{Department of Mathematics and Computer Science,
University of Basel, Spiegelgasse 1, 4051 Basel, Switzerland,
(marcus.grote@unibas.ch)}
\and Simon R. J. Michel\thanks{Mathematisches Institut, Universit\"{a}t Bern,
Sidlerstrasse 5, 3012 Bern, Switzerland, (simon.michel@unibe.ch)}
\and Stefan A. Sauter\thanks{Institute for Mathematics, University of Zurich,
Winterthurerstrasse 190, 8057 Zurich, Switzerland, (stas@math.uzh.ch)} }
\maketitle
\date{}

\begin{abstract}
Adaptivity and local mesh refinement are crucial for the efficient numerical simulation 
of wave phenomena in complex geometry. Local mesh refinement, however, can impose
a tiny time-step across the entire computational domain when using explicit
time integration.
By taking smaller time-steps yet only inside locally refined regions, local time-stepping methods
overcome the stringent CFL stability restriction imposed on the global time-step by
a small fraction of the elements without sacrificing explicitness.  
In \cite{Grote_Mitkova}, a leapfrog based local time-stepping method
was proposed for the inhomogeneous wave equation, which
applies standard leapfrog time-marching with a
smaller time-step inside the refined region. 
Here, to remove potential instability at certain time-steps, a stabilized version is proposed which leads
to optimal $L^2$-error estimates under a CFL condition independent of the coarse-to-fine mesh
ratio. Moreover, a weighted transition is introduced 
to restore optimal $H^1$-convergence when the source is nonzero
across the coarse-to-fine mesh interface. Numerical experiments corroborate the theoretical error estimates and illustrate
the usefulness of these improvements.
\end{abstract}

\section{Introduction}
For the simulation of wave phenomena,
finite element methods, be it conforming or discontinuous, 
easily accommodate complex geometry or material interfaces. Local
mesh refinement, however, seriously hampers the efficiency
of explicit time integrators, such as the popular leapfrog (LF) method, 
when a small fraction of elements impose a tiny time-step for stability
across the entire mesh. By taking smaller time-steps, but only inside 
those refined regions, while taking larger time-steps elsewhere, local
time-stepping methods (LTS) permit to overcome that overly
stringent CFL stability constraint on the global time-step
without sacrificing the explicitness or parallel efficiency
of explicit time integrators. 

Following the classical IMEX (implicit-explicit) approach, first locally implicit methods that combine 
the (explicit) second-order LF (or St\o rmer-Verlet) method in the coarse part with the (implicit) Crank-Nicolson
method in the fine part of the mesh
were proposed and further analyzed in  \cite{Piperno,DFFL10,Ver2011,DLM13,HS16,HS19}.
In \cite{ChabassierImperiale2015}, fourth-order energy-preserving IMEX schemes were derived for the wave equation, when  the computational domain is divided
 into a coarse and fine region using Lagrange multipliers along the interface.

Even earlier, an explicit local time-stepping (LTS) method was proposed for the wave equation, 
which also employs the explicit LF method in the fine part of the mesh, yet with a smaller 
step-size~\cite{ColFouJol3,ColFouJol1}; error estimates were derived
in~\cite{ColFouJol2, JolyRodriguez} in the one-dimensional case. Stability of 
this second-order method, however, is guaranteed by enforcing energy conservation 
which requires at every time-step the solution of a linear
system for the shared unknowns at the coarse/fine intersection; hence, the overall scheme is not fully explicit. 

Fully explicit second-order LTS methods based on St\o rmer-Verlet in time and DG discretizations in space were proposed for Maxwell's equations in \cite{Piperno,MPFC08}. Starting from the popular LF method, fully explicit
local time-stepping methods of arbitrary order were proposed
for the homogeneous wave equation in \cite{DG09} and
for the inhomogeneous wave equation with nonzero forcing in \cite{Grote_Mitkova}.
Inside the locally refined region, they apply standard leapfrog time-marching with a
smaller time-step, which also implies evaluating the inhomogeneous right-hand side
at all intermediate times. 
Hence, when combined with a mass-lumped conforming \cite{COHEN,mass_lumping_2d} or discontinuous Galerkin FE discretization \cite{GSS06} in space, the resulting LF-LTS method is truly explicit, inherently parallel and also conserves a discrete energy in the absence of forcing; it was successfully applied to 3D seismic wave propagation \cite{MZKM13}. A multilevel version was later proposed \cite{DG15} and achieved high parallel efficiency \cite{Rietmann2017}.

In \cite{grote_sauter_1}, optimal $L^2$-error estimates for the LF-LTS method from \cite{DG09} were
derived for $H^1$-conforming FEM discretizations (with mass-lumping), though under a CFL condition
where the global time-step $\Delta t$ in fact depended on the smallest elements in the mesh.
For the convergence analysis, the inner loop over the $p$ local LF steps of size $\Delta t / p$ was rewritten
in terms of a single global time-step $\Delta t$, which involves Chebyshev polynomials of the first kind.
In general one cannot improve upon the stability constraint on 
$\Delta t$, as the LF-LTS method may become unstable at certain values of $\Delta t$. 

By replacing those Chebyshev polynomials by 
``damped'' Chebyshev polynomials with damping parameter $\nu >0$, which slightly modifies the LF-LTS 
scheme inside the refined region, a stabilized LF-LTS($\nu$) version was proposed in \cite{grote2021stabilized} for the homogeneous wave equation, and also independently in \cite{CarleHochbruck2}, which led to optimal $L^2$-error estimates yet under a CFL condition
independent of the coarse-to-fine aspect ratio. The same damped Chebyshev polynomials were also previously used
to stabilize leapfrog-Chebyshev (LFC, discrete Gautschi-type) methods \cite{GilbertJoly,CarleHochbruchCheby}, where the ``non-stiff'' (or slowly varying), possibly nonlinear source term is evaluated only at every global time-step throughout the entire computational domain. In fact, 
the error analysis for semilinear ODE's  in \cite{CarleHochbruck2} and for DG discretizations of linear wave equations
 in \cite{CarleHochbruck3} both apply to the inhomogeneous case, too, when the source is evaluated only once during each global time-step, akin to the approach in \cite{CarleHochbruchCheby}.

{\bf Main contributions.} 
The leapfrog based local time-stepping (LF-LTS) method introduced 
in \cite{Grote_Mitkova} for the inhomogeneous wave equation can become unstable
for certain time-steps. To remove those
potential instabilities, we propose a stabilized version which leads
to optimal $L^2$-error estimates under a CFL condition independent of the coarse-to-fine mesh
ratio. Thus we effectively
introduce stabilization into the original method from \cite{Grote_Mitkova}, which, in contrast to the ``split-LFC'' method from \cite{CarleHochbruck3}, corresponds to standard LF time-marching, yet with a smaller time-step, inside the refined region. Generalizing the proof of $L^2$-convergence for the homogeneous wave equation  \cite{grote2021stabilized} to the inhomogeneous case is not 
straightforward due to the source evaluations at all intermediate local time-steps. In fact it requires the development of new 
analytical techniques involving Chebyshev polynomials of the second kind
to rewrite and analyze the effect of those local time-steps and thus
finally prove the stability and optimal consistency order of the method.

In \cite{ChabassierImperiale2021} it was shown that the expected optimal $H^1$-convergence rate of the 
LF-LTS method \cite{Grote_Mitkova} for the inhomogeneous wave equation 
is reduced by $h^{1/2}$ when the source is nonzero across the coarse-to-fine mesh interface. To overcome this
difficulty and thereby restore optimal $H^1$-convergence, we introduce a new weighted transition which relies on a
discrete mesh-based distance function.
 
{\bf Outline.}
The rest of our article is structured as follows. In Section \ref{SecGalDisc}, we first recall the $H^1$-conforming Galerkin formulation of the inhomogeneous wave equation together with needed notation as the mesh is split into a coarse and a locally refined part. 
In Section \ref{SecLTSalgo}, this leads to the stabilized LF-LTS$(\nu)$ algorithm for the inhomogenous wave equation, which also evaluates
the source term inside the locally refined region at all intermediate time-steps.  
In Section \ref{SecConvAna}, we prove optimal $L^2$-convergence rates under a CFL condition independent of the fine-to-coarse mesh ratio. To accomodate for the
intermediate source evaluations, our error analysis requires both Chebyshev polynomials of the first and the second kind.
Still our analysis also includes the somewhat simpler situation when the source is evaluated at global time-steps only as in \cite{CarleHochbruck3}. To achieve optimal $H^1$-convergence regardless of the source, 
we introduce in Section \ref{SecH1conv} the new weighted transition which relies on a
discrete mesh-based distance function. In Section \ref{SecNumEx}, we present numerical experiments which illustrate the 
optimal $L^2$-convergence rates expected from theory. 
In doing so, 
we also compare the LF-LTS method with the alternative ``split-LFC'' approach from \cite{CarleHochbruck3}, where the
source inside the locally refined region is only evaluated once at each global time-step.
While both LTS methods, either with or without intermediate local
source evaluations, are second-order convergent in $\Delta t$, our approach is more accurate (but also slightly more expensive) than the LTS method from \cite{CarleHochbruck3} when
the source strongly varies locally. 
Finally, we revisit the original example from  \cite{ChabassierImperiale2021} 
and show that the new weighted transition restores the expected optimal $H^1$-convergence rate even when the source
is nonzero across the coarse-to-fine mesh interface.

\section{Stabilized LF-LTS Method\label{SecGalDisc}}

In Sections \ref{sec:waveeq} and \ref{sec:FEM}, we adopt and recall the setting for finite element spaces as in \cite{grote2021stabilized}.

\subsection{The Wave Equation} \label{sec:waveeq}

Let $\Omega\subset\mathbb{R}^{d}$ be a bounded Lipschitz domain with boundary
$\Gamma=\partial\Omega$. For $1\leq q\leq\infty$, let $L^{q}\left(
\Omega\right)  $ be the standard Lebesgue space with norm $\left\Vert
\cdot\right\Vert _{L^{q}\left(  \Omega\right)  }$. For $q=2$, the scalar
product in $L^{2}\left(  \Omega\right)  $ is denoted by $\left(  \cdot
,\cdot\right)  $ and the norm by $\left\Vert \cdot\right\Vert :=\left\Vert
\cdot\right\Vert _{L^{2}\left(  \Omega\right)  }$. For $k\in\mathbb{N}_{0}$
and $1\leq q\leq\infty$, let $W^{k,q}\left(  \Omega\right)  $ denote the
classical Sobolev spaces equipped with the norm $\left\Vert \cdot\right\Vert
_{W^{k,q}\left(  \Omega\right)  }$. For $q=2$, these spaces are Hilbert spaces
and denoted by $H^{k}\left(  \Omega\right)  :=W^{k,2}\left(  \Omega\right)  $
with scalar product $\left(  \cdot,\cdot\right)  _{H^{k}\left(  \Omega\right)
}$ and norm $\left\Vert \cdot\right\Vert _{H^{k}\left(  \Omega\right)
}:=\left(  \cdot,\cdot\right)  _{H^{k}\left(  \Omega\right)  }^{1/2}%
=\left\Vert \cdot\right\Vert _{W^{k,2}\left(  \Omega\right)  }$. 

We assume that the boundary is split into a Dirichlet and a Neumann part:
$\Gamma=\Gamma_{D}\cup\Gamma_{N}$, $\Gamma_{D}\cap\Gamma_{N}=\emptyset$ with
the convention that $\Gamma_{D}$ is (relatively) closed and has positive
surface measure. On $\Gamma_{D}$ homogeneous Dirichlet boundary conditions are
imposed and Neumann boundary conditions on $\Gamma_{N}$. The space of Sobolev
functions in $H^{1}\left(  \Omega\right)  $ which vanish on the Dirichlet part
of the boundary $\partial\Omega$ is denoted by $H_{D}^{1}\left(
\Omega\right)  :=\left\{  w\in H^{1}\left(  \Omega\right)  :\left.
w\right\vert _{\Gamma_{D}}=0\right\}  $. For $\Gamma_{D}=\Gamma$, we use the
shorthand $H_{0}^{1}\left(  \Omega\right)  $. Throughout this paper we
restrict ourselves to function spaces over the field of real numbers.

Let $V\subset H^{1}\left(  \Omega\right)  $ denote a closed subspace of
$H^{1}\left(  \Omega\right)  $, such as $V=H^{1}\left(  \Omega\right)  $ or
$V=H_{0}^{1}\left(  \Omega\right)  $, and $a:V\times V\rightarrow\mathbb{R}$
denote a bilinear form, which is symmetric, continuous, and coercive:%
\begin{subequations}
\label{wellposed}
\end{subequations}%
\begin{equation}
a\left(  u,v\right)  =a\left(  v,u\right)  \qquad\forall u,v\in V\tag{%
\ref{wellposed}%
a}\label{wellposeda}%
\end{equation}
and%
\begin{equation}
\left\vert a\left(  u,v\right)  \right\vert \leq C_{\operatorname*{cont}%
}\left\Vert u\right\Vert _{H^{1}\left(  \Omega\right)  }\left\Vert
v\right\Vert _{H^{1}\left(  \Omega\right)  }\qquad\forall u,v\in V\tag{%
\ref{wellposed}%
b}\label{wellposedb}%
\end{equation}
and%
\begin{equation}
a\left(  u,u\right)  \geq c_{\operatorname*{coer}}\left\Vert u\right\Vert
_{H^{1}\left(  \Omega\right)  }^{2}\qquad\forall u\in V.\tag{%
\ref{wellposed}%
c}\label{wellposedc}%
\end{equation}
For given $u_{0}\in V,v_{0}\in L^{2}\left(  \Omega\right)  $ and $F:\left[
0,T\right]  \rightarrow V^{\prime}$, we consider the wave equation: Find
$u:\left[  0,T\right]  \rightarrow V$ such that%
\begin{equation}
\left(  \ddot{u},w\right)  +a\left(  u,w\right)  =F\left(  w\right)
\quad\forall w\in V,t>0\label{waveeq}%
\end{equation}
with initial conditions%
\begin{equation}
u\left(  0\right)  =u_{0}\quad\text{and\quad}\dot{u}\left(  0\right)
=v_{0}.\label{waveeqic}%
\end{equation}
It is well known that (\ref{waveeq})--(\ref{waveeqic}) is well-posed for
sufficiently regular $u_{0}$, $v_{0}$ and $F$ \cite{LionsMagenesI}. In fact,
the weak solution $u$ can be shown to be continuous in time, that is, $u\in
C^{0}(0,T;V),\dot{u}\in C^{0}(0,T;L^{2}\left(  \Omega\right)  )$ -- see
\cite[Chapter III, Theorems 8.1 and 8.2]{LionsMagenesI} for details -- which
implies that the initial conditions (\ref{waveeqic}) are well defined.
Moreover, we always assume that there exists a function $f:\left[  0,T\right]
\rightarrow L^{2}\left(  \Omega\right)  $ such that%
\[
F\left(  t\right)  \left(  w\right)  =\left(  f\left(  t\right)  ,w\right)
\qquad\forall w\in V\quad\forall t\in\left[  0,T\right]  .
\]
For the stability and convergence analysis we will impose further smoothness
assumptions on $f$.

\begin{example}
\label{Exmodel problem}The classical second-order wave equation in strong form
is given by%
\begin{equation}%
\begin{split}
\partial_{t}^{2}u-\nabla\cdot(c^{2}\nabla u)  &  =f\qquad\;\,\mbox{in }\Omega
\times(0,T),\\
u  &  =0\qquad\;\,\mbox{on }\Gamma_{D}\times(0,T),\\
\frac{\partial u}{\partial\nu}  &  =0\qquad\;\,\mbox{on }\Gamma_{N}%
\times(0,T),\\
u|_{t=0}  &  =u_{0}\qquad\mbox{in }\Omega,\\
u_{t}|_{t=0}  &  =v_{0}\qquad\mbox{in }\Omega,\\
&
\end{split}
\label{model problem}%
\end{equation}
where the velocity field $c(x)$ satisfies $0<c_{\min}\leq c(x)\leq c_{\max}$.
In this case, we have $V:=H_{D}^{1}\left(  \Omega\right)  :=\left\{  w\in
H^{1}\left(  \Omega\right)  :\left.  w\right\vert _{\Gamma_{D}}=0\right\}  $;
the bilinear form is given by $a\left(  u,v\right)  :=\left(  c^{2}\nabla
u,\nabla u\right)  $ and the right-hand side by $F\left(  w\right)  =\left(
f,w\right)  $ for all $w\in V$.
\end{example}

\subsection{Galerkin Finite Element Spatial Discretization} \label{sec:FEM}

We begin by briefly introducing our notation for $hp$ finite element
spaces.\ Let $\mathcal{T}$ be a simplicial finite element mesh for $\Omega$
which is conforming in the sense that there are no \textit{hanging nodes.} The
local and global mesh width is%
\[
h_{\tau}:=\operatorname*{diam}\tau,\quad h:=\max_{\tau\in\mathcal{T}}h_{\tau
}.
\]
As a convention the triangles in $\mathcal{T}$ are closed sets. Without loss
of generality, we may assume that there is fixed $h_{0}>0$ such that%
\begin{equation}
h\leq h_{0}. \label{hsmallerh0}%
\end{equation}
The finite element space of continuous, piecewise polynomials of degree
$m\geq1$ is denoted by
\[
S_{\mathcal{T}}^{m}:=\left\{  u\in C^{0}\left(  \Omega\right)  \mid\forall
\tau\in\mathcal{T}\quad\left.  u\right\vert _{\tau}\in\mathbb{P}_{m}\right\}
,
\]
where $\mathbb{P}_{m}$ is the space of polynomials of maximal total degree
$m$. The standard local nodal basis is denoted by $b_{m,z}$, $z\in
\Sigma_{\mathcal{T}}^{m}$, where $\Sigma_{\mathcal{T}}^{m}$ is the full set of
nodal points of order $m$.
The finite element space with incorporated Dirichlet boundary conditions is
given by $S:=S_{\mathcal{T}}^{m}\cap V$. The set of nodal points for $S$ is
the subset $\Sigma:=\left\{  z\in\Sigma_{\mathcal{T}}^{m}\mid z\in
\Gamma_{\operatorname*{D}}\right\}  $ and the basis in $S$ is denoted by
$b_{z}:=b_{m,z}$, $z\in\Sigma$. Furthermore, we use the notion of a
mass-lumped scalar product%
\[
\left(  u,v\right)  _{\mathcal{T}}:=\sum_{z\in\Sigma}d_{z}u\left(  z\right)
v\left(  z\right)  ,\quad\forall u,v\in S
\]
for certain coefficients $d_{z}>0$ and induced norm $\Vert\cdot\Vert
_{\mathcal{T}}:=\left(  \cdot,\cdot\right)  _{\mathcal{T}}^{1/2}$ which is
assumed to be equivalent to the $L^{2}\left(  \Omega\right)  $-norm: there
exists a constant $C_{\operatorname*{eq}}>0$ such that%
\begin{equation}
C_{\operatorname*{eq}}^{-1}\left\Vert u\right\Vert _{\mathcal{T}}%
\leq\left\Vert u\right\Vert \leq C_{\operatorname*{eq}}\left\Vert u\right\Vert
_{\mathcal{T}}\quad\forall u\in S\text{.} \label{ceqCeq}%
\end{equation}
For more details we refer, e.g., to \cite[Definition 2.3]{grote2021stabilized}.

We also associate to the discrete bilinear form the linear operator
$A^{S}:S\rightarrow S$ defined by%
\begin{equation}
\left(  A^{S}u,v\right)  _{\mathcal{T}}=a\left(  u,v\right)  \qquad\forall
u,v\in S \label{defAS}%
\end{equation}
and introduce the discrete right-hand side function $f_{S}:[0,T]\rightarrow S$
defined by
\[
\left(  f_{S}(t),v\right)  _{\mathcal{T}}=F(t)(v)\qquad\forall\,v\in
S\quad\forall\,t\in\lbrack0,T].
\]

With the definitions and notation at hand, the semi-discrete wave equation
(with mass-lumping) is given by: find $u_{S}:\left[  0,T\right]  \rightarrow
S$ such that%
\begin{subequations}
\label{spacedisc}
\end{subequations}%
\begin{equation}
\left(  \ddot{u}_{S},v\right)  _{\mathcal{T}}+a\left(  u_{S},v\right)
=\left(  f_{S}(t),v\right)  _{\mathcal{T}}\quad\forall v\in S,t>0\tag{%
\ref{spacedisc}%
a}\label{spacedisca}%
\end{equation}
with initial conditions%
\begin{equation}
u_{S}\left(  0\right)  =r_{S}u_{0}\quad\text{and\quad}\dot{u}_{S}\left(
0\right)  =r_{S}v_{0},\tag{%
\ref{spacedisc}%
b}\label{spacediscb}%
\end{equation}
where $r_{S}u_{0}$ is a finite element approximation of $u_{0}$ (see, e.g.,
\cite[(2.8)]{grote2021stabilized}).

As in \cite{grote2021stabilized}, we consider finite element meshes
$\mathcal{T}$ which consist of a quasi-uniform part $\mathcal{T}%
_{\operatorname*{c}}\subset\mathcal{T}$ and allow for local refinements so
that $\mathcal{T}_{\operatorname*{f}}:=\mathcal{T}\backslash\mathcal{T}%
_{\operatorname*{c}}$ might contain much smaller triangles than $\mathcal{T}%
_{\operatorname{c}}$. In this way, we assume that the mesh size of
$\mathcal{T}_{\operatorname{c}}$ satisfies $h_{\operatorname{c}}=h$.
Similarly, we define coarse and locally refined regions of the
domain\footnote{By $\operatorname*{int}\left(  M\right)  $ we denote the open
interior of a subset $M\subset\mathbb{R}^{d}$.} as
\begin{equation}
\Omega_{\operatorname*{c}}:=\operatorname*{int}\left(  {\textstyle\bigcup
\nolimits_{\tau\in\mathcal{T}_{\operatorname*{c}}}}\tau\right)  ,\quad
\Omega_{\operatorname*{f}}:=\operatorname*{int}\left(  {\textstyle\bigcup
\nolimits_{\tau\in\mathcal{T}_{\operatorname*{f}}}}\tau\right)  ,\quad
\Omega_{\operatorname*{f}}^{+1}:=\operatorname*{int}\left(  {\textstyle\bigcup
\nolimits_{\substack{\tau\in\mathcal{T}\\\tau\cap\overline{\Omega
_{\operatorname*{f}}}\neq\emptyset}}}\,\tau\right)  .\label{defsubdomains}%
\end{equation}
Hence, $\Omega_{\operatorname*{f}}^{+1}$ contains $\Omega_{\operatorname*{f}}$
and all elements directly adjacent to it. We note that $\Omega
_{\operatorname*{f}}$ and $\Omega_{\operatorname*{c}}$ are disjoint, while
their union covers all of $\Omega$ except for the coarse/fine interface, and
that $\Omega_{\operatorname*{c}}\cup\Omega_{\operatorname*{f}}^{+1}=\Omega$
since the interface between $\Omega_{\operatorname*{f}}$ and $\Omega
_{\operatorname*{c}}$ is contained in $\Omega_{\operatorname*{f}}^{+1}$.

We also split the degrees of freedom associated with the fine or coarse parts
of the mesh, respectively, as
\[
\Sigma_{\operatorname*{f}}:=\Sigma\cap\overline{\Omega_{\operatorname*{f}}%
}\quad\text{and\quad}\Sigma_{\operatorname*{c}}:=\Sigma\backslash
\Sigma_{\operatorname*{f}},
\]
and introduce the corresponding FE subspaces
\[
S_{\operatorname*{c}}:=\operatorname{span}\left\{  b_{z}:z\in\Sigma
_{\operatorname*{c}}\right\}  \quad\mbox{and}\quad S_{\operatorname*{f}%
}:=\operatorname{span}\left\{  b_{z}:z\in\Sigma_{\operatorname*{f}}\right\}  ,
\]
which gives rise to a direct sum decomposition: for every $u\in S$ there exist
unique functions $u_{\operatorname*{c}}\in S_{\operatorname*{c}}$ and
$u_{\operatorname*{f}}\in S_{\operatorname*{f}}$ such that
\begin{equation}
u=u_{\operatorname*{c}}+u_{\operatorname*{f}}.\label{uadddecomp}%
\end{equation}
Hence, we can define the projections $\Pi_{\operatorname*{c}}^{S}:S\rightarrow
S_{\operatorname*{c}}$ and $\Pi_{\operatorname*{f}}^{S}:S\rightarrow
S_{\operatorname*{f}}$ by
\begin{equation}
\Pi_{\operatorname*{c}}^{S}u:=u_{\operatorname*{c}}\quad\text{and}\quad
\Pi_{\operatorname*{f}}^{S}u:=u_{\operatorname*{f}}.\label{defPifPic}%
\end{equation}
In fact, the decomposition (\ref{uadddecomp}) is orthogonal with respect to
the $\left(  \cdot,\cdot\right)  _{\mathcal{T}}$ scalar product \cite[Lemma
3.2]{grote2021stabilized}, which is essential to prove sharp bounds for the
eigenvalues of our discrete bilinearform \cite[Theorem 3.7]%
{grote2021stabilized}.

\subsection{Stabilized LF-LTS($\nu$) Algorithm} \label{SecLTSalgo}

The original leapfrog (LF) based local time-stepping (LF-LTS$(0)$) methods for the
numerical solution of the second-order wave equations \eqref{waveeq} were
proposed for homogeneous right-hand sides in \cite{DG09} and for inhomogeneous
right-hand sides in \cite{Grote_Mitkova}.
Optimal convergence rates for the
LF-LTS$(0)$ method from \cite{DG09} with $p$ local time steps were derived for a
conforming FEM discretization, albeit under a CFL condition where $\Delta t$
in fact depends on the smallest elements in the mesh \cite{grote_sauter_1}. To
prove optimal $L^{2}$ convergence rates under a CFL condition independent of
$p$, a stabilized algorithm LF-LTS$(\nu)$ was introduced recently in
\cite{grote2021stabilized}, and also independently in \cite{CarleHochbruck2}. 
Here, $\nu \geq 0$ denotes a small stabilization parameter; typically, we set $\nu=0.01$.
Stability and convergence for homogeneous right-hand sides were derived in
\cite{grote2021stabilized} for a CFL condition independent of the
coarse/fine mesh ratio. 

We shall now extend this stabilized LF-LTS$(\nu)$ method to inhomogeneous
right-hand sides so that for $\nu=0$ the method coincides with the original
LF-LTS$(0)$ method \cite[Sect. 4.1]{Grote_Mitkova} for the inhomogeneous wave equation.
To do so, let the (global) time-step $\Delta t=T/N$ and $u_{S}^{\left(
n\right)  }$ denote the FE approximation at time $t_{n}=n\Delta t$. Given the
numerical solution at times $t_{n-1}$ and $t_{n}$, the LF-LTS method then
computes the numerical solution of \eqref{spacedisc} at $t_{n+1}$ by using a
smaller time-step $\Delta\tau=\Delta t/p$ inside the regions of local
refinement; here, $p\geq2$ denotes the \textquotedblleft
coarse\textquotedblright\ to \textquotedblleft fine\textquotedblright%
\ time-step ratio.

Inside the fine region, the source $f$ is also evaluated over the time
interval $[t_n - \Delta t, t_n + \Delta t]$ at the intermediate
times $t_{n+k/p}=\left(  n+k/p\right)  \Delta t$ and we let
\[
f_{S,k}^{(n)}:=f_{S}\left(  t_{n}+\frac{k}{p}\Delta t\right)  ,\quad-p\leq
k\leq p.
\]

The cost of those extra evaluations inside the locally refined region is typically negligible compared to the
overall solution process. Still, a locally strongly varying source typically calls for such a 
a higher sampling of the right-hand side -- see Section \ref{SecL2convNumEx}.

Given the fixed stabilization parameter $0\leq
\nu\leq1/2$, we now define the following constants determined by 
stabilized versions \cite{HV03} of Chebyshev
polynomials of the first and second kind \cite{Rivlin}, denoted by $T_{n}$ and
$U_{n}$, respectively:
\begin{equation}
\delta_{p,\nu}:=1+\frac{\nu}{p^{2}},\qquad\omega_{p,\nu}:=2\,\frac
{T_{p}^{\prime}\left(  \delta_{p,\nu}\right)  }{T_{p}\left(  \delta_{p,\nu
}\right)  },\label{defdeltaomega}%
\end{equation}
and%
\begin{equation}
\beta_{p,\nu}^{(k,\ell)}:=\frac{T_{k+\ell}\left(  \delta_{p,\nu}\right)
}{T_{k+1}\left(  \delta_{p,\nu}\right)  }\quad\text{and}\quad\gamma_{p,\nu
}^{(k)}:=\frac{\left(  p-k\right)  \beta_{p,\nu}^{\left(  k,p-k\right)  }%
}{U_{p-1-k}\left(  \delta_{p,\nu}\right)  }\label{Def_beta}%
\end{equation}
for $k=0,1,\ldots,p-1$ and $\ell=-1,0,\ldots,p-k$; formally we set $T_{-1}=0$.
Later the special case
\begin{equation}
\beta_{p,\nu}^{(0,0)}=\frac{\gamma_{p,\nu}^{\left(  0\right)  }}{2p^{2}}%
\omega_{p,\nu}=\frac{1}{T_{1}\left(  \delta_{p,\nu}\right)  }=\frac{1}%
{\delta_{p,\nu}}\label{betpnuerel}%
\end{equation}
will be employed which follows from (\ref{Def_beta}) and the formula
$T_{n}^{\prime}=nU_{n-1}$ (see \cite[18.9.21]{NIST:DLMF}).

The stabilized LF-LTS algorithm for the inhomogeneous wave equation
which computes the new solution $u_{S}^{(n+1)}$ from the current and previous
values $u_{S}^{(n)}$ and $u_{S}^{(n-1)}$, is
now given by:

\begin{algorithm}
[LF-LTS$(\nu)$ Galerkin FE Algorithm]\label{AlgStab} Let $n\geq1$.

\begin{enumerate}
\item Given $u_{S}^{(n-1)}$, $u_{S}^{(n)}$, compute $w_{S}^{(n)}$ as%
\begin{equation}
w_{S}^{(n)}=\Pi_{\operatorname*{c}}^{S}f_{S,0}^{(n)}-A^{S}\Pi
_{\operatorname*{c}}^{S}u_{S}^{(n)}. \label{Defwsn}%
\end{equation}

\item Compute
\begin{equation}
z_{S,1}^{(n)}=u_{S}^{\left(  n\right)  }+\frac{1}{2}\left(  \frac{\Delta t}%
{p}\right)  ^{2}\left(  \frac{2p^{2}}{\omega_{p,\nu}}\,\beta_{p,\nu}%
^{(0,0)}\left(  w_{S}^{(n)}-A^{S}\Pi_{\operatorname*{f}}^{S}u_{S}%
^{(n)}\right)  +\gamma_{p,\nu}^{(0)}\Pi_{\operatorname*{f}}^{S}f_{S,0}%
^{(n)}\right)  . \label{AlgStep2}%
\end{equation}

\item For $k=1,\ldots,p-1,$ compute
\begin{align*}
z_{S,k+1}^{(n)}  &  =\left(  1+\beta_{p,\nu}^{(k,-1)}\right)  z_{S,k}%
^{(n)}-\beta_{p,\nu}^{(k,-1)}z_{S,k-1}^{(n)}\\
&  \quad+\left(  \frac{\Delta t}{p}\right)  ^{2}\left(  \frac{2p^{2}}%
{\omega_{p,\nu}}\,\beta_{p,\nu}^{(k,0)}\left(  w_{S}^{(n)}-A^{S}%
\Pi_{\operatorname*{f}}^{S}z_{S,k}^{(n)}\right)  +\gamma_{p,\nu}^{(k)}\frac
{1}{2}\Pi_{\operatorname*{f}}^{S}\left(  f_{S,k}^{(n)}+f_{S,-k}^{(n)}\right)
\right)  .
\end{align*}

\item Compute
\[
u_{S}^{(n+1)}=-u_{S}^{(n-1)}+2z_{S,p}^{(n)}.
\]

\end{enumerate}
\end{algorithm}

The $\operatorname*{LF-LTS}(\nu)$ method is completed by specifying the initial
values
\begin{align}
u_{S}^{\left(  0\right)  }  & =r_{S}u_{0},\nonumber\\
u_{S}^{\left(  1\right)  }  & =r_{S}u_{0}+\Delta t\,r_{S}v_{0}+\frac{\Delta
t^{2}}{2}\left(  f_{S}^{(0)}-A^{S}u_{0}\right)  .\tag{%
\ref{eq:StabLFLTSrhs}%
b}\label{eq:StabLFLTSrhs_init}%
\end{align}
Here, the term $A^{S}u_{0}$ in the second equation could be replaced by
$A^{S,p,\nu}u_{0}$ thereby allowing for local time-stepping already during the
very first time-step. In that case, the analysis below also applies but
requires a minor change as explained in \cite[Rem. 2.5]{grote2021stabilized}.
This modification neither affects the stability nor the convergence rate of
the overall LF-LTS$(\nu)$ scheme.

\begin{remark}
\label{rem:LFCLTScomparison}
The above LF-LTS algorithm evaluates the source $f$ inside the refined region
at all intermediate times $t_{n+k/p}=t_n + (k/p)  \Delta t$, $k = 1-p, \dots, p-1$.
Hence for $\nu=0$, it coincides with the original LF-LTS method
for the inhomogeneous wave equation \cite[Sect. 4.1]{Grote_Mitkova}, since
\[
\delta_{p,0}=1,\quad\omega_{p,0}=2p^{2},\quad\beta_{p,0}^{(k,\ell)}%
=1,\quad\gamma_{p,0}^{(k)}=\frac{p-k}{U_{p-1-k}(1)}=1
\]
holds for all $k=0,1,\ldots,p-1$, $\ell=-1,0,\ldots,p-k$. In particular, for $\nu = 0$ it 
corresponds to a standard LF discretization of the inhomogeneous
wave equation with the smaller time-step $\Delta \tau = \Delta t/ p$. In contrast, the
``split-LFC'' local time integration method in \cite{CarleHochbruck3} omits those intermediate evaluations of
$f$ inside the fine region, akin to the previous LFC method \cite{CarleHochbruchCheby}. That second variant
from \cite{CarleHochbruck3} without intermediate source evaluations is easily obtained by replacing
in the above algorithm  \eqref{Defwsn} with
\[
w_S^{(n)} = f_{S,0}^{(n)}-A^{S}\Pi
_{\operatorname*{c}}^{S}u_{S}^{(n)}
\]
and setting $\gamma_{p,\nu}^{(k)} = 0$ for $k=0,1,\ldots,p-1 .$
Clearly if $f=0$, both algorithms coincide. 
\end{remark}

\section{Convergence Theory in the $L^{2}$-norm}
\label{SecConvAna}

\subsection{Two-Step Formulation via Chebyshev Polynomials}

For the analysis, we will first rewrite the above LF-LTS$(\nu)$ algorithm in a
two-step leapfrog-like formulation as it was done in
\cite{grote2021stabilized} for homogeneous right-hand sides. For this, we
define the polynomials%
\begin{align}
P_{p,\nu,k}\left(  x\right)   &  :=2\left(  1-\frac{T_{k}\left(  \delta
_{p,\nu}-\frac{x}{\omega_{p,\nu}}\right)  }{T_{k}\left(  \delta_{p,\nu
}\right)  }\right)  ,\label{defPp}\\
P_{p,\nu,k}^{\Delta t}\left(  x\right)   &  :=\frac{P_{p,\nu,k}\left(  \Delta
t^{2}x\right)  }{\Delta t^{2}x},\label{defPpDeltat}\\
Q_{p,\nu,r,k}^{\Delta t}\left(  x\right)   &  :=\frac{p-k}{p^{2}}\frac
{T_{p}\left(  \delta_{p,\nu}\right)  }{T_{k}\left(  \delta_{p,\nu}\right)
}\frac{U_{k-1-r}\left(  \delta_{p,\nu}-\frac{\Delta t^{2}x}{\omega_{p,\nu}%
}\right)  }{U_{p-1-r}\left(  \delta_{p,\nu}\right)  },\label{eq:defQ_pnurk}%
\end{align}
where formally we set $U_{-1}=0$ and write short $P_{p,\nu}:=P_{p,\nu,p}$,
$P_{p,\nu}^{\Delta t}:=P_{p,\nu,p}^{\Delta t}$, and $Q_{p,\nu,r}^{\Delta
t}:=Q_{p,\nu,r,p}^{\Delta t}$ in what follows.

\begin{remark}
\label{RemSimplePropPQ}The univariate polynomials $P_{p,\nu,k}^{\Delta t}$,
$Q_{p,\nu,r,k}^{\Delta t}$ will appear in an explicit representation formula
for the functions $z_{S,k}^{(n)}$ at the local time steps. Some simple
properties are listed below:%
\[
P_{p,\nu,k}^{\Delta t}\in\mathbb{P}_{k-1}\qquad Q_{p,\nu,r,k}^{\Delta t}%
\in\mathbb{P}_{k-1-r}%
\]
and the first two terms in a Taylor series around $x=0$ are given by%
\begin{align}
P_{p,\nu,k}^{\Delta t}\left(  x\right)   &  =\frac{2}{\omega_{p,\nu}}%
\frac{T_{k}^{\prime}\left(  \delta_{p,\nu}\right)  }{T_{k}\left(
\delta_{p,\nu}\right)  }-x\frac{\Delta t^{2}}{\omega_{p,\nu}^{2}}\frac
{T_{k}^{\prime\prime}\left(  \delta_{p,\nu}\right)  }{T_{k}\left(
\delta_{p,\nu}\right)  }+\ldots\label{firsttermexp}\\
Q_{p,\nu,r,k}^{\Delta t}\left(  x\right)   &  =\frac{p-k}{p^{2}}\frac
{T_{p}\left(  \delta_{p,\nu}\right)  }{T_{k}\left(  \delta_{p,\nu}\right)
}\left(  1-x\frac{\Delta t^{2}}{\omega_{p,\nu}}\frac{U_{k-1-r}^{\prime}\left(
\delta_{p,\nu}\right)  }{U_{p-1-r}\left(  \delta_{p,\nu}\right)  }%
+\ldots\right)  .\label{firsttermexp2}%
\end{align}
Note that the expression for $Q_{p,\nu,r,k}^{\Delta t}\left(  x\right)  $ is
valid only for $r\leq k-1$ while for $k=r$ we have $Q_{p,\nu,r,r+1}^{\Delta
t}=0$. Since $\delta_{p,\nu}$ lies outside the interval $\left]  -1,1\right[
$ the denominators in these expressions cannot be zero. For the special case
$k=1$ it is straightforward to verify%
\begin{equation}
P_{p,\nu,1}^{\Delta t}(x)=Q_{p,\nu,0,1}^{\Delta t}(x)=\frac{\gamma_{p,\nu
}^{(0)}}{p^{2}}.\label{Ppnue1}%
\end{equation}
\end{remark}

The following recurrence relations for $P_{p,\nu,k}$ follow easily from the
recurrence relations for Chebyshev polynomials as already stated in \cite[Lem.
6.1]{CarleHochbruchCheby} and \cite[Lem. B.1]{grote2021stabilized}.
\begin{lemma}
Let $r=0,1,\ldots,p-1$. Then, the polynomials $P_{p,\nu,k}$ defined in
(\ref{defPp}) satisfy the recurrence relation%
\begin{align}
P_{p,\nu,0}\left(  x\right)   &  =0,\nonumber\\
P_{p,\nu,1}\left(  x\right)   &  =\frac{2x}{\delta_{p,\nu}\omega_{p,\nu}%
},\label{Ppnuex}\\
P_{p,\nu,k+1}\left(  x\right)   &  =2\beta_{p,\nu}^{\left(  k,0\right)
}\left(  \delta_{p,\nu}-\frac{x}{\omega_{p,\nu}}\right)  P_{p,\nu,k}\left(
x\right)  -\beta_{p,\nu}^{\left(  k,-1\right)  }P_{p,\nu,k-1}\left(  x\right)
+\frac{4}{\omega_{p,\nu}}\beta_{p,\nu}^{\left(  k,0\right)  }x,\nonumber
\end{align}
while $P_{p,\nu,k}^{\Delta t}$ satisfies%
\begin{align}
P_{p,\nu,0}^{\Delta t}(x)  &  =0,\nonumber\\
P_{p,\nu,1}^{\Delta t}(x)  &  =\frac{2}{\delta_{p,\nu}\omega_{p,\nu}%
},\label{Ppnuedtx}\\
P_{p,\nu,k+1}^{\Delta t}(x)  &  =2\beta_{p,\nu}^{\left(  k,0\right)  }\left(
\delta_{p,\nu}-\frac{\Delta t^{2}x}{\omega_{p,\nu}}\right)  P_{p,\nu
,k}^{\Delta t}\left(  x\right)  -\beta_{p,\nu}^{\left(  k,-1\right)  }%
P_{p,\nu,k-1}^{\Delta t}\left(  x\right)  +\frac{4}{\omega_{p,\nu}}%
\beta_{p,\nu}^{\left(  k,0\right)  }. \label{Pdeltarek}%
\end{align}

\end{lemma}

\begin{lemma}
\label{lem:recursion_Q_pnurk} Let $r=0,1,\ldots,p-1$. Then, the polynomials
$Q_{p,\nu,r,k}^{\Delta t}$ defined in \eqref{eq:defQ_pnurk} satisfy the
recurrence relation%
\begin{align}
Q_{p,\nu,r,r}^{\Delta t}(x)  &  =0, &  & \label{eq:recursion_Q_pnurk_0}\\
Q_{p,\nu,r,r+1}^{\Delta t}(x)  &  =\frac{\gamma_{p,\nu}^{(r)}}{p^{2}}, &  &
\label{eq:recursion_Q_pnurk_1}\\
Q_{p,\nu,r,k+1}^{\Delta t}(x)  &  =2\beta_{p,\nu}^{(k,0)}\,\left(
\delta_{p,\nu}-\frac{\Delta t^{2}}{\omega_{p,\nu}}x\right)  Q_{p,\nu
,r,k}^{\Delta t}(x) &  & \nonumber\\
&  \quad-\beta_{p,\nu}^{(k,-1)}Q_{p,\nu,r,k-1}^{\Delta t}(x), &  &  r+1\leq
k\leq p-1. \label{eq:recursion_Q_pnurk_k}%
\end{align}

\end{lemma}%

\proof
The first equation \eqref{eq:recursion_Q_pnurk_0} holds by
\eqref{eq:defQ_pnurk} and \eqref{eq:recursion_Q_pnurk_1} follows directly
since $U_{0}\equiv1$. Hence, let $r+1\leq k\leq p-1$ and insert definitions
(\ref{Def_beta}), (\ref{eq:defQ_pnurk}) into the right-hand side of
(\ref{eq:recursion_Q_pnurk_k}) to obtain%
\begin{align*}
&  \left(  1+\beta_{p,\nu}^{(k,-1)}-\left(  \frac{\Delta t}{p}\right)
^{2}\frac{2p^{2}}{\omega_{p,\nu}}\beta_{p,\nu}^{(k,0)}\,x\right)
Q_{p,\nu,r,k}^{\Delta t}(x)-\beta_{p,\nu}^{(k,-1)}Q_{p,\nu,r,k-1}^{\Delta
t}\left(  x\right)  \\
&  =\left(  \frac{T_{k+1}\left(  \delta_{p,\nu}\right)  +T_{k-1}\left(
\delta_{p,\nu}\right)  -2\frac{\Delta t^{2}x}{\omega_{p,\nu}}T_{k}\left(
\delta_{p,\nu}\right)  }{T_{k}\left(  \delta_{p,\nu}\right)  }U_{k-1-r}\left(
\delta_{p,\nu}-\frac{\Delta t^{2}x}{\omega_{p,\nu}}\right)  \right.  \\
&  \qquad\left.  -U_{k-2-r}\left(  \delta_{p,\nu}-\frac{\Delta t^{2}x}%
{\omega_{p,\nu}}\right)  \right)  \frac{\gamma_{p,\nu}^{(r)}}{p^{2}}%
\frac{T_{r+1}\left(  \delta_{p,\nu}\right)  }{T_{k+1}\left(  \delta_{p,\nu
}\right)  }.
\end{align*}
Using the standard recursions for Chebyshev polynomials of first and second
kind \cite[Table 18.9.1]{NIST:DLMF}, respectively, we conclude
\begin{align*}
&  \left(  1+\beta_{p,\nu}^{(k,-1)}-\left(  \frac{\Delta t}{p}\right)
^{2}\frac{2p^{2}}{\omega_{p,\nu}}\beta_{p,\nu}^{(k,0)}\,x\right)
Q_{p,\nu,r,k}^{\Delta t}(x)-\beta_{p,\nu}^{(k,-1)}Q_{p,\nu,r,k-1}^{\Delta
t}(x)\\
&  =\frac{\gamma_{p,\nu}^{(r)}}{p^{2}}\frac{T_{r+1}\left(  \delta_{p,\nu
}\right)  }{T_{k+1}\left(  \delta_{p,\nu}\right)  }\left(  2\left(
\delta_{p,\nu}-\frac{\Delta t^{2}x}{\omega_{p,\nu}}\right)  U_{k-1-r}\left(
\delta_{p,\nu}-\frac{\Delta t^{2}x}{\omega_{p,\nu}}\right)  -U_{k-2-r}\left(
\delta_{p,\nu}-\frac{\Delta t^{2}x}{\omega_{p,\nu}}\right)  \right)  \\
&  =\frac{\gamma_{p,\nu}^{(r)}}{p^{2}}\frac{T_{r+1}\left(  \delta_{p,\nu
}\right)  }{T_{k+1}\left(  \delta_{p,\nu}\right)  }U_{k-r}\left(
\delta_{p,\nu}-\frac{\Delta t^{2}x}{\omega_{p,\nu}}\right)  \\
&  =Q_{p,\nu,r,k+1}^{\Delta t}(x).
\end{align*}%
\endproof

\begin{lemma}
For $0\leq k\leq p$, the functions $z_{S,k}^{(n)}$ defined in Algorithm
\ref{AlgStab} can be written in the form
\begin{equation}
z_{S,k}^{(n)}=\mathfrak{L}_{S,k}^{(n)}u_{S}^{(n)}+\frac{\Delta t^{2}}%
{2}\left(  \frac{2}{\omega_{p,\nu}}\frac{T_{k}^{\prime}\left(  \delta_{p,\nu
}\right)  }{T_{k}\left(  \delta_{p,\nu}\right)  }\Pi_{\operatorname*{c}}%
^{S}f_{S,0}^{(n)}+\sum_{r=-(k-1)}^{k-1}Q_{p,\nu,|r|,k}^{\Delta t}\left(
A^{S}\Pi_{\operatorname*{f}}^{S}\right)  \Pi_{\operatorname*{f}}^{S}%
f_{S,r}^{(n)}\right)  \label{eq:zn_Sk_explicit}%
\end{equation}
for%
\[
\mathfrak{L}_{S,k}^{(n)}=\left(  I^{S}-\frac{\Delta t^{2}}{2}P_{p,\nu
,k}^{\Delta t}\left(  A^{S}\Pi_{\operatorname*{f}}^{S}\right)  A^{S}\right)  .
\]
\label{lem:zSkn_explrepr}
\end{lemma}%

\proof
In \cite[Lemma B.2, (B.2)]{grote2021stabilized}, relation
(\ref{eq:zn_Sk_explicit}) is shown for $f_{S}=0$. In the inhomogeneous case
here, it follows by the same arguments that
\begin{equation}
z_{S,k}^{(n)}=\mathfrak{L}_{S,k}^{(n)}u_{S}^{(n)}+\frac{\Delta t^{2}}%
{2}\mathfrak{R}_{S,k}^{(n)},\label{eq:z_Sk_leftright}%
\end{equation}
where the term $\mathfrak{R}_{S,k}^{(n)}$ only depends on the right-hand side
evaluations $\left(  f_{S,r}^{(n)}\right)
_{r=-\left(  k-1\right)  }^{k-1}$.

Hence, it remains to prove that%
\begin{equation}
\mathfrak{R}_{S,k}^{(n)}=\frac{2}{\omega_{p,\nu}}\frac{T_{k}^{\prime}\left(
\delta_{p,\nu}\right)  }{T_{k}\left(  \delta_{p,\nu}\right)  }\Pi
_{\operatorname*{c}}^{S}f_{S,0}^{(n)}+\sum_{r=-(k-1)}^{k-1}Q_{p,\nu
,|r|,k}^{\Delta t}\left(  A^{S}\Pi_{\operatorname*{f}}^{S}\right)
\Pi_{\operatorname*{f}}^{S}f_{S,r}^{\left(  n\right)  }.\label{eq:PropTS}%
\end{equation}
From $\Pi_{\operatorname{f}}^{S}\Pi_{\operatorname{c}}^{S}f=0$ (see
(\ref{defPifPic})) it follows:%
\[
P_{p,\nu,k}^{\Delta t}\left(  A^{S}\Pi_{\operatorname*{f}}^{S}\right)
\Pi_{\operatorname*{c}}^{S}f_{S,0}^{(n)}=P_{p,\nu,k}^{\Delta t}\left(
0\right)  \Pi_{\operatorname*{c}}^{S}f_{S,0}^{(n)}\overset
{\text{(\ref{firsttermexp})}}{=}\frac{2}{\omega_{p,\nu}}\frac{T_{k}^{\prime
}\left(  \delta_{p,\nu}\right)  }{T_{k}\left(  \delta_{p,\nu}\right)  }%
\Pi_{\operatorname*{c}}^{S}f_{S,0}^{(n)}%
\]
so that
\begin{equation}
\mathfrak{R}_{S,k}^{(n)}=P_{p,\nu,k}^{\Delta t}\left(  A^{S}\Pi
_{\operatorname*{f}}^{S}\right)  \Pi_{\operatorname*{c}}^{S}f_{S,0}^{(n)}%
+\sum_{r=-(k-1)}^{k-1}Q_{p,\nu,|r|,k}^{\Delta t}\left(  A^{S}\Pi
_{\operatorname*{f}}^{S}\right)  \Pi_{\operatorname*{f}}^{S}f_{S,r}^{\left(
n\right)  }\label{eq:PropTS2}%
\end{equation}
implies (\ref{eq:PropTS}). 

We prove (\ref{eq:PropTS2}) by induction. For $k=0$, this is trivial, and for
$k=1$, we use the definition of $z_{S,1}^{\left(  n\right)  }$ in Algorithm
\ref{AlgStab}, insert (\ref{eq:recursion_Q_pnurk_1}) for $r=0$ (see also
Remark \ref{RemSimplePropPQ}), employ $w_{S}^{(n)}=\Pi_{\operatorname*{c}}%
^{S}f_{S,0}^{(n)}-A^{S}\Pi_{\operatorname*{c}}^{S}u_{S}^{(n)}$ (see
(\ref{Defwsn})) and (\ref{eq:z_Sk_leftright}) to obtain%
\begin{align*}
z_{S,1}^{(n)} &  =u_{S}^{\left(  n\right)  }+\frac{1}{2}\left(  \frac{\Delta
t}{p}\right)  ^{2}\left(  \frac{2p^{2}}{\omega_{p,\nu}}\,\beta_{p,\nu}%
^{(0,0)}\left(  w_{S}^{(n)}-A^{S}\Pi_{\operatorname*{f}}^{S}u_{S}%
^{(n)}\right)  +\gamma_{p,\nu}^{(0)}\Pi_{\operatorname*{f}}^{S}f_{S,0}%
^{(n)}\right)  \\
&  =\mathfrak{L}_{S,1}^{(n)}u_{S}^{(n)}+\frac{1}{2}\left(  \frac{\Delta t}%
{p}\right)  ^{2}\left(  \frac{2p^{2}}{\omega_{p,\nu}}\,\beta_{p,\nu}%
^{(0,0)}\Pi_{\operatorname*{c}}^{S}f_{S,0}^{(n)}+\gamma_{p,\nu}^{(0)}%
\Pi_{\operatorname*{f}}^{S}f_{S,0}^{(n)}\right)  \\
&  \overset{\text{(\ref{eq:recursion_Q_pnurk_1})}}{=}\mathfrak{L}_{S,1}%
^{(n)}u_{S}^{(n)}+\frac{\Delta t^{2}}{2}\left(  \frac{2}{\omega_{p,\nu}}%
\beta_{p,\nu}^{(0,0)}\Pi_{\operatorname*{c}}^{S}f_{S,0}^{(n)}+Q_{p,\nu
,0,1}^{\Delta t}\left(  A^{S}\Pi_{\operatorname*{f}}^{S}\right)
\Pi_{\operatorname*{f}}^{S}f_{S,0}\right)  .
\end{align*}
From $\beta_{p,\nu}^{\left(  0,0\right)  }=\frac{T_{1}^{\prime}\left(
\delta_{p,\nu}\right)  }{T_{1}\left(  \delta_{p,\nu}\right)  }$ (cf.
(\ref{betpnuerel})) it follows that the second summand is $\frac{\Delta t^{2}%
}{2}\mathfrak{R}_{S,1}^{(n)}$.

For $k>1$, Step 3 in Algorithm \ref{AlgStab} in combination with the induction
hypothesis, the definition of $w_{S}^{\left(  n\right)  }$, and the fact that
the linear part with respect $u_{S}^{\left(  n\right)  }$ is known from
(\ref{eq:z_Sk_leftright}) leads to%
\begin{align*}
z_{S,k+1}^{(n)} &  =\left(  1+\beta_{p,\nu}^{(k,-1)}\right)  z_{S,k}%
^{(n)}-\beta_{p,\nu}^{(k,-1)}z_{S,k-1}^{(n)}\\
&  \quad+\left(  \frac{\Delta t}{p}\right)  ^{2}\left(  \frac{2p^{2}}%
{\omega_{p,\nu}}\,\beta_{p,\nu}^{(k,0)}\left(  w_{S}^{(n)}-A^{S}%
\Pi_{\operatorname*{f}}^{S}z_{S,k}^{(n)}\right)  +\gamma_{p,\nu}^{(k)}\frac
{1}{2}\Pi_{\operatorname*{f}}^{S}\left(  f_{S,k}^{(n)}+f_{S,-k}^{(n)}\right)
\right)  \\
&  =\mathfrak{L}_{S,k+1}^{(n)}u_{S}^{(n)}+\frac{\Delta t^{2}}{2}\rho
\end{align*}
for%
\begin{align*}
\rho:=\,
&2\beta_{p,\nu}^{(k,0)}\left(  \delta_{p,\nu}-\frac{\Delta t^{2}%
}{\omega_{p,\nu}}A^{S}\Pi_{\operatorname*{f}}^{S}\right)  \mathfrak{R}%
_{S,k}^{(n)}-\beta_{p,\nu}^{(k,-1)}\mathfrak{R}_{S,k-1}^{(n)} \\
&\quad +4\frac
{\,\beta_{p,\nu}^{(k,0)}}{\omega_{p,\nu}}\Pi_{\operatorname*{c}}^{S}%
f_{S,0}^{(n)}+\frac{\gamma_{p,\nu}^{(k)}}{p^{2}}\Pi_{\operatorname*{f}}%
^{S}\left(  f_{S,k}^{(n)}+f_{S,-k}^{(n)}\right)  ,
\end{align*}
where we used
\[
\left(  1+\beta_{p,\nu}^{(k,-1)}\right)  /\left(  2\beta_{p,\nu}%
^{(k,0)}\right)  =\frac{T_{k+1}\left(  \delta_{p,\nu}\right)  +T_{k-1}\left(
\delta_{p,\nu}\right)  }{2T_{k}\left(  \delta_{p,\nu}\right)  }\overset
{\text{\cite[Table 18.9.1]{NIST:DLMF}}}{=}\delta_{p,\nu}.
\]
It remains to prove that $\rho=\mathfrak{R}_{S,k+1}^{(n)}$. We use the
induction hypotheses for $\mathfrak{R}_{S,k}^{(n)}$ as in (\ref{eq:PropTS}),
(\ref{eq:recursion_Q_pnurk_1}) and obtain%
\begin{align*}
\rho &  :=\left.  \left(  2\beta_{p,\nu}^{(k,0)}\left(  \delta_{p,\nu}%
-\frac{\Delta t^{2}}{\omega_{p,\nu}}x\right)  P_{p,\nu,k}^{\Delta t}\left(
x\right)  -\beta_{p,\nu}^{(k,-1)}P_{p,\nu,k-1}^{\Delta t}\left(  x\right)
+4\frac{\,\beta_{p,\nu}^{(k,0)}}{\omega_{p,\nu}}\right)  \right\vert
_{x=A^{S}\Pi_{\operatorname*{f}}^{S}}\Pi_{\operatorname*{c}}^{S}f_{S,0}%
^{(n)}\\
&  +\sum_{r=-(k-1)}^{k-1}\left.  \left(  2\beta_{p,\nu}^{(k,0)}\left(
\delta_{p,\nu}-\frac{\Delta t^{2}}{\omega_{p,\nu}}x\right)  Q_{p,\nu
,|r|,k}^{\Delta t}\left(  x\right)  -\beta_{p,\nu}^{(k,-1)}Q_{p,\nu
,|r|,k-1}^{\Delta t}\left(  x\right)  \right)  \right\vert _{x=A^{S}%
\Pi_{\operatorname*{f}}^{S}}\Pi_{\operatorname*{f}}^{S}f_{S,r}^{\left(
n\right)  }\\
&  +Q_{p,\nu,k,k+1}^{\Delta t}\left(  A^{S}\Pi_{\operatorname*{f}}^{S}\right)
\Pi_{\operatorname*{f}}^{S}\left(  f_{S,k}^{(n)}+f_{S,-k}^{(n)}\right)
\end{align*}
where we used $Q_{p,\nu,k-1,k-1}^{\Delta t}=0$. The first row can be
simplified by using the recurrence for $P_{p,\nu,k+1}^{\Delta t}$ (see
(\ref{Ppnuex})). For the second one we employ the recursion
(\ref{eq:recursion_Q_pnurk_k}) for $Q_{p,\nu,|r|,k}^{\Delta t}$ so that the
second row equals%
\begin{align*}
& \sum_{r=-(k-1)}^{k-1}\left.  Q_{p,\nu,r,k+1}^{\Delta t}\left(  x\right)
\right\vert _{x=A^{S}\Pi_{\operatorname*{f}}^{S}}\Pi_{\operatorname*{f}}%
^{S}f_{S,r}^{\left(  n\right)  }+Q_{p,\nu,k,k+1}^{\Delta t}\left(  A^{S}%
\Pi_{\operatorname*{f}}^{S}\right)  \Pi_{\operatorname*{f}}^{S}\left(
f_{S,k}^{(n)}+f_{S,-k}^{(n)}\right)  \\
& \qquad=\sum_{r=-k}^{k}\left.  Q_{p,\nu,r,k+1}^{\Delta t}\left(  x\right)
\right\vert _{x=A^{S}\Pi_{\operatorname*{f}}^{S}}\Pi_{\operatorname*{f}}%
^{S}f_{S,r}^{\left(  n\right)  }.
\end{align*}
In this way we have proved (\ref{eq:PropTS2}) and, in turn, (\ref{eq:PropTS}).%
\endproof

\begin{corollary}
\label{CorLTS2step}
Let $a^{p,\nu}:S\times S\rightarrow\mathbb{R}$ denote the discrete bilinear
form
\begin{equation}
a^{p,\nu}\left(  u,v\right)  :=a\left(  P_{p,\nu}^{\Delta t}\left(
\Pi_{\operatorname*{f}}^{S}A^{S}\right)  u,v\right)  \qquad\forall u,v\in
S\label{defapv}%
\end{equation}
with associated (linear) operator $A^{S,p,\nu}:S\rightarrow S$%
\begin{align}
A^{S,p,\nu} &  :=\left(  A^{S}\right)  ^{1/2}P_{p,\nu}^{\Delta t}\left(
\left(  A^{S}\right)  ^{1/2}\Pi_{\operatorname*{f}}^{S}\left(  A^{S}\right)
^{1/2}\right)  \left(  A^{S}\right)  ^{1/2}\label{defASp}\\
&  =P_{p,\nu}^{\Delta t}\left(  A^{S}\Pi_{\operatorname*{f}}^{S}\right)
A^{S}=A^{S}P_{p,\nu}^{\Delta t}\left(  \Pi_{\operatorname*{f}}^{S}%
A^{S}\right)  \nonumber
\end{align}
and let $R_{S}:[0,T]\rightarrow S$ by defined by%
\begin{equation}
R_{S}\left(  t\right)  :=P_{p,\nu}^{\Delta t}\left(  A^{S}\Pi
_{\operatorname*{f}}^{S}\right)  \Pi_{\operatorname*{c}}^{S}f_{S}%
(t)+\sum_{r=-(p-1)}^{p-1}Q_{p,\nu,|r|}^{\Delta t}\left(  A^{S}\Pi
_{\operatorname*{f}}^{S}\right)  \Pi_{\operatorname*{f}}^{S}f_{S}\left(
t+\frac{r\Delta t}{p}\right)  .\label{eq:SumRhs}%
\end{equation}
Then, for $n\geq1$, the $\operatorname*{LF-LTS}(\nu)$ algorithm is equivalent to%
\begin{subequations}
\label{eq:StabLFLTSrhs}
\end{subequations}%
\begin{equation}
\left(  u_{S}^{(n+1)}-2\,u_{S}^{(n)}+u_{S}^{(n-1)},w\right)  _{\mathcal{T}%
}=\Delta t^{2}\left[  \left(  R_{S}^{(n)},w\right)  _{\mathcal{T}}-a^{p,\nu
}\left(  u_{S}^{(n)},w\right)  \right]  \quad\forall\,w\in S,\tag{%
\ref{eq:StabLFLTSrhs}%
a}\label{eq:StabLFLTSrhs_n}%
\end{equation}
where $R_{S}^{(n)}:=R_{S}\left(  t_{n}\right)  $.
\end{corollary}%

\proof
The assertion follows directly from inserting (\ref{eq:zn_Sk_explicit}) for
$k=p$ into Step 4 of Algorithm \ref{AlgStab}.%
\endproof

\subsection{Main Convergence Results}

The time-steps $\Delta t$ and $\Delta\tau=\Delta t/p$, $p\in\mathbb{N}$, used
in $\mathcal{T}_{\operatorname*{c}}$ and $\mathcal{T}_{\operatorname*{f}}$,
respectively, are each determined by the smallest triangle in either part of
the mesh. Hence, the \textquotedblleft coarse\textquotedblright%
-to-\textquotedblleft fine\textquotedblright\ time-step ratio, $p$, must
satisfy
\[
\dfrac{\min\left\{  h_{\tau}:\tau\in\mathcal{T}_{\operatorname*{c}}\right\}
}{h_{\tau}}\leq p\qquad\forall\tau\in\mathcal{T}.
\]
Given the quasi-uniformity constant for the \emph{coarse} mesh $\mathcal{T}%
_{\operatorname*{c}}$,
\[
C_{\operatorname*{qu},\operatorname*{c}}:=\frac{h_{\operatorname*{c}}}%
{\min\left\{  h_{\tau}:\tau\in\mathcal{T}_{\operatorname*{c}}\right\}  },
\]
the ratio between the largest and smallest elements of the mesh is thus
bounded by
\begin{equation}
\dfrac{h_{\operatorname*{c}}}{h_{\tau}}\leq C_{\operatorname*{qu}%
,\operatorname*{c}}p\qquad\forall\tau\in\mathcal{T}.\label{hfphc}%
\end{equation}
In \cite[Section 3.1]{grote2021stabilized}, theoretical properties for the
bilinear form $a^{p,\nu}$ are derived, which will be used in the error
analysis below. We impose the same CFL-type stability condition as in
\cite[(3.8)]{grote2021stabilized} on the ratio $\Delta t/h_{\operatorname*{c}%
}$:
\begin{equation}
\left(  3+\frac{C_{\operatorname*{cont}}}{c_{\operatorname*{coer}}}\right)
C_{\operatorname*{cont}}C_{\operatorname*{inv}}^{2}\left(  \frac{\Delta
t}{h_{\operatorname*{c}}}\right)  ^{2}\leq\frac{\nu}{\nu+1},\label{CFLtotal}%
\end{equation}
where $C_{\operatorname*{inv}}$ is a constant independent of $p$ such that the
inverse inequalities holds (cf. \cite[(3.3), (3.4)]{grote2021stabilized}):%
\begin{equation}
\left.
\begin{array}
[c]{c}%
\left\Vert u\right\Vert _{H^{1}\left(  \Omega_{\operatorname{c}}\right)  }\leq
C_{\operatorname*{inv}}h_{\operatorname*{c}}^{-1}\left\Vert u\right\Vert
_{\mathcal{T}}\\
\left\Vert u\right\Vert _{H^{1}\left(  \Omega\right)  }\leq
C_{\operatorname*{inv}}\dfrac{p}{h_{\operatorname*{c}}}\left\Vert u\right\Vert
_{\mathcal{T}}%
\end{array}
\right\}  \quad\forall u\in S.\label{DefCinv}%
\end{equation}

For a discussion and a possible simplification of the CFL condition
(\ref{CFLtotal}) we refer to \cite[(3.8)]{grote2021stabilized}.

{To derive a priori error estimates for the LF-LTS$(\nu)$ Galerkin FE solution of
\eqref{eq:StabLFLTSrhs}, we first introduce%
\[
v_{S}^{\left(  n+1/2\right)  }:=\frac{u_{S}^{\left(  n+1\right)  }%
-u_{S}^{\left(  n\right)  }}{\Delta t},
\]
and rewrite \eqref{eq:StabLFLTSrhs} as a one-step method
\begin{equation}
\begin{aligned} \left( v_{S}^{\left( n+1/2\right) },q\right)_{\mathcal{T}} & = \left( v_{S}^{\left( n-1/2\right) },q\right)_{\mathcal{T}} -\Delta t\, a^{p,\nu}\left( u_{S}^{\left( n\right) },q\right) + \Delta t \left( R_S^{(n)},q\right)_{\mathcal{T}} \; & &\forall q\in S,\\ \left( u_{S}^{\left( n+1\right) },r\right)_{\mathcal{T}} & = \left( u_{S}^{\left( n\right) },r\right)_{\mathcal{T}} + \Delta t\left( v_{S}^{\left( n+1/2\right) },r\right)_{\mathcal{T}} \; & &\forall r\in S,\\ \left( u_{S}^{\left( 0\right) },w\right)_{\mathcal{T}} & = \left( u_{0},w\right)_{\mathcal{T}} & &\forall w\in S,\\ \left( v_{S}^{\left( 1/2\right) },w\right)_{\mathcal{T}} & =\left( v_{0},w\right)_{\mathcal{T}} + \frac{\Delta t}{2} \left( \left( f_S^{(0)},w\right)_{\mathcal{T}} - a\left( u_{0},w\right) \right) \;& &\forall w\in S. \end{aligned} \label{eq1}%
\end{equation}
}

The first two equations in (\ref{eq1}) correspond to the one-step iteration
\[
\left[
\begin{array}
[c]{c}%
v_{S}^{(n+1/2)}\\
u_{S}^{(n+1)}%
\end{array}
\right]  =\mathfrak{B}\left[
\begin{array}
[c]{c}%
v_{S}^{(n-1/2)}\\
u_{S}^{(n)}%
\end{array}
\right]  +\Delta t%
\begin{bmatrix}
R_{S}^{(n)}\\
\Delta t\,R_{S}^{(n)}%
\end{bmatrix}
\quad n\geq1
\]
with
\[
\mathfrak{B}:=\left[
\begin{array}
[c]{cc}%
I_{S} & -\Delta tA^{S,p,\nu}\\
\Delta t\,I_{S} & I_{S}-\Delta t^{2}A^{S,p,\nu}%
\end{array}
\right]
\]
and $A^{S,p,\nu}$ as in (\ref{defASp}).

Next, we denote the error by
\[
\mathbf{e}^{\left(  n+1\right)  }:=%
\begin{bmatrix}
e_{v}^{\left(  n+1/2\right)  }\\
e_{u}^{\left(  n+1\right)  }%
\end{bmatrix}
:=%
\begin{bmatrix}
v\left(  t_{n+1/2}\right)  -v_{S}^{\left(  n+1/2\right)  }\\
u\left(  t_{n+1}\right)  -u_{S}^{\left(  n+1\right)  }%
\end{bmatrix}
,
\]
where $u$ is the solution of \eqref{waveeq}-\eqref{waveeqic} and $v$ the
solution of the corresponding first-order formulation: {Find
$u,v:[0,T]\rightarrow V$ such that%
\begin{align*}
\left(  \dot{v},w\right)  +a\left(  u,w\right)   &  = \left(f\left(t\right),w\right) & \quad\forall\,w &
\in V,\quad t>0,\\
\left(  v,w\right)   &  =\left(  \dot{u},w\right)   & \quad\forall\,w &  \in
V,\quad t>0
\end{align*}
with initial conditions $u(0)=u_{0}$ and $v(0)=v_{0}$.}

We shall split the error into a semi-discrete and a fully discrete
contribution. To do so, we introduce the first-order formulation of the
semi-discrete problem (\ref{spacedisc}). Find $u_{S},v_{S}:\left[  0,T\right]
\rightarrow S$ such that%
\[%
\begin{array}
[c]{cl}%
\left.
\begin{array}
[c]{rcl}%
\left(  \dot{v}_{S},w\right)  _{\mathcal{T}}+a\left(  u_{S},w\right)   & = &
\left(f_S\left(t\right),w\right)_{\mathcal{T}}\\
\left(  v_{S},w\right)  _{\mathcal{T}} & = & \left(  \dot{u}_{S},w\right)
_{\mathcal{T}}%
\end{array}
\right\}   & \forall w\in S,\quad t>0,\\%
\begin{array}
[c]{rcl}%
\qquad u_{S}\left(  0\right)   & = & r_{S}u_{0}\\
\qquad v_{S}\left(  0\right)   & = & r_{S}v_{0}.
\end{array}
&
\end{array}
\]
Hence, we may write
\begin{equation}
\mathbf{e}^{\left(  n+1\right)  }=\mathbf{e}_{S}^{\left(  n+1\right)
}+\mathbf{e}_{S,\Delta t}^{\left(  n+1\right)  }\label{defspliterror}%
\end{equation}
with%
\begin{align*}
\mathbf{e}_{S}^{\left(  n+1\right)  } &  :=\left[
\begin{array}
[c]{c}%
e_{v,S}^{\left(  n+1/2\right)  }\\
e_{u,S}^{\left(  n+1\right)  }%
\end{array}
\right]  :=\left[
\begin{array}
[c]{c}%
v\left(  t_{n+1/2}\right)  -v_{S}\left(  t_{n+1/2}\right)  \\
u\left(  t_{n+1}\right)  -u_{S}\left(  t_{n+1}\right)
\end{array}
\right]  ,\\
\mathbf{e}_{S,\Delta t}^{\left(  n+1\right)  } &  :=\left[
\begin{array}
[c]{c}%
e_{v,S,\Delta t}^{\left(  n+1/2\right)  }\\
e_{u,S,\Delta t}^{\left(  n+1\right)  }%
\end{array}
\right]  :=\left[
\begin{array}
[c]{c}%
v_{S}\left(  t_{n+1/2}\right)  -v_{S}^{\left(  n+1/2\right)  }\\
u_{S}\left(  t_{n+1}\right)  -u_{S}^{\left(  n+1\right)  }%
\end{array}
\right]  .
\end{align*}

Following similar arguments as in \cite[\S 3.2]{grote_sauter_1}, we derive a
recurrence relation for $\mathbf{e}_{S,\Delta t}^{(n+1)}$, which can be solved
and eventually yields the explicit error representation
\begin{align}
\left[
\begin{array}
[c]{c}%
e_{v,S,\Delta t}^{\left(  n+1/2\right)  }\\
e_{u,S,\Delta t}^{\left(  n+1\right)  }%
\end{array}
\right]   &  =\mathfrak{B}^{n}\left[
\begin{array}
[c]{c}%
e_{v,S,\Delta t}^{\left(  1/2\right)  }\\
e_{u,S,\Delta t}^{\left(  1\right)  }%
\end{array}
\right]  +\Delta t\sum_{\ell=1}^{n-1}\mathfrak{B}^{\ell}\left(
\mbox{\boldmath$ \sigma$}^{\left(  n-\ell\right)  }%
-\mbox{\boldmath$ \sigma$}^{\left(  n-\ell+1\right)  }\right)
\label{errorreprfu}\\
&  \quad+\Delta t\mbox{\boldmath$ \sigma$}^{\left(  n\right)  }-\Delta
t\mathfrak{B}^{n}\mbox{\boldmath$ \sigma$}^{\left(  1\right)  },\nonumber
\end{align}
where
\begin{equation}
\mbox{\boldmath$ \sigma$}^{\left(  n\right)  }=\frac{1}{\Delta t}%
\begin{bmatrix}
-\dfrac{u_{S}\left(  t_{n+1}\right)  -u_{S}\left(  t_{n}\right)  }{\Delta
t}+v_{S}\left(  t_{n+1/2}\right)  \\
u_{S}\left(  t_{n}\right)  +\left(  A^{S,p,\nu}\right)  ^{-1}\left(
\dfrac{v_{S}\left(  t_{n+1/2}\right)  -v_{S}\left(  t_{n-1/2}\right)  }{\Delta
t}-R_{S}^{(n)}\right)
\end{bmatrix}
\label{defboldsigmab}%
\end{equation}
with $R_S^{(n)}$ defined as in Corollary \ref{CorLTS2step}. 

Before proving the main convergence theorem we provide some approximation estimates.

\begin{lemma}
\label{lem:errorestRn} Let $R_{S}$ be defined as in (\ref{eq:SumRhs}) and let
$\theta\in\left[  t_{n},t_{n+1}\right]  $, $\nu > 0$. Then, there exist constants $C_{1}$,
$C_{2}$ independent of $\Delta t$, $h$, $p$, and $f_{S}$ such that%
\begin{equation}
\left\Vert \left(  A^{S,p,\nu}\right)  ^{-1}\partial_{t}^{3}R_{S}\left(
\theta\right)  \right\Vert _{\mathcal{T}}\leq C_{1}\max_{t\in\lbrack
t_{n-1},t_{n+2}]}{\left\Vert \partial_{t}^{3}f_{S}\left(  t\right)
\right\Vert _{\mathcal{T}}}\label{lem:errorestRnEst1}%
\end{equation}
and
\begin{equation}
\left\Vert \left(  A^{S,p,\nu}\right)  ^{-1}\partial_{t}^{\ell}\left(
f_{S}^{(n)}-R_{S}^{(n)}\right)  \right\Vert _{\mathcal{T}}\leq C_{2}\Delta
t^{2}\max_{k=0,2}\left\{  \max_{t\in\lbrack t_{n},t_{n+1}]}{\left\Vert
\partial_{t}^{k+\ell}f_{S}\left(  t\right)  \right\Vert _{\mathcal{T}}%
}\right\}  \label{lem:errorestRnEst2}%
\end{equation}
for $\ell=0,1$.
\end{lemma}%

\proof
For an operator $K:S\rightarrow S$ we introduce the operator norm%
\[
\left\Vert K\right\Vert _{S\leftarrow S}:=\sup_{w\in S\backslash\left\{
0\right\}  }\frac{\left\Vert Kw\right\Vert _{\mathcal{T}}}{\left\Vert
w\right\Vert _{\mathcal{T}}}.
\]
We start from the definition%
\begin{align}
\left(  A^{S,p,\nu}\right)  ^{-1}\partial_{t}^{\ell}R_{S}\left(  t\right)    &
=\left(  A^{S,p,\nu}\right)  ^{-1}P_{p,\nu}^{\Delta t}\left(  A^{S}%
\Pi_{\operatorname*{f}}^{S}\right)  \Pi_{\operatorname*{c}}^{S}\partial
_{t}^{\ell}f_{S}\left(  t\right)  \label{Rsl}\\
& +\left(  A^{S,p,\nu}\right)  ^{-1}\sum_{r=-(p-1)}^{p-1}Q_{p,\nu,|r|}^{\Delta
t}\left(  A^{S}\Pi_{\operatorname*{f}}^{S}\right)  \Pi_{\operatorname*{f}}%
^{S}\partial_{t}^{\ell}f_{S}\left(  t+\frac{r\Delta t}{p}\right)  .\nonumber
\end{align}
Since $\Pi_{\operatorname{f}}^{S}\Pi_{\operatorname{c}}^{S}f_{S}=0$ we obtain
from (\ref{firsttermexp}) and (\ref{defdeltaomega}) for the first term in
(\ref{Rsl})%
\begin{equation}
P_{p,\nu}^{\Delta t}\left(  A^{S}\Pi_{\operatorname*{f}}^{S}\right)
\Pi_{\operatorname*{c}}^{S}\partial_{t}^{\ell}f_{S}\left(  t\right)
=\Pi_{\operatorname*{c}}^{S}\partial_{t}^{\ell}f_{S}\left(  t\right)
.\label{Ppnue2}%
\end{equation}
and hence%
\[
\left\Vert P_{p,\nu}^{\Delta t}\left(  A^{S}\Pi_{\operatorname*{f}}%
^{S}\right)  \Pi_{\operatorname*{c}}^{S}\partial_{t}^{\ell}f_{S}\left(
t\right)  \right\Vert _{\mathcal{T}}\leq\left\Vert \partial_{t}^{\ell}%
f_{S}\left(  t\right)  \right\Vert _{\mathcal{T}}%
\]
From \cite[Theorem 3.7]{grote2021stabilized} it follows%
\[
\left\Vert \left(  A^{S,p,\nu}\right)  ^{-1}v\right\Vert _{\mathcal{T}}%
\leq\frac{1}{c_{\nu}}\left\Vert v\right\Vert _{\mathcal{T}},
\]
with%
\begin{equation}
c_{\nu}=\frac{1}{2}\min\left\{  \frac{c_{\operatorname*{coer}}}%
{C_{\operatorname*{eq}}^{2}},\frac{2\nu}{\left(  \nu+1\right)  \Delta t^{2}%
}\right\}  \label{Defcnue}%
\end{equation}
and, in turn,%
\begin{equation}
\left\Vert \left(  A^{S,p,\nu}\right)  ^{-1}\partial_{t}^{3}P_{p,\nu}^{\Delta
t}\left(  A^{S}\Pi_{\operatorname*{f}}^{S}\right)  \Pi_{\operatorname*{c}}%
^{S}\partial_{t}^{\ell}f_{S}\left(  t\right)  \right\Vert _{\mathcal{T}}%
\leq\frac{1}{c_{\nu}}\left\Vert \partial_{t}^{\ell}f_{S}\left(  t\right)
\right\Vert _{\mathcal{T}}.
\label{eq:estPpnue2cnue}
\end{equation}
Next we consider the second summand in (\ref{Rsl}). The decomposition operator
$E^{S}:S\rightarrow S_{\operatorname{c}}\times S_{\operatorname{f}}$ is
defined by $E^{S}u=\left(  \Pi_{\operatorname{c}}^{S}u,\Pi_{\operatorname{f}%
}^{S}u\right)  $ (cf. (\ref{defPifPic})) so that its left-inverse
$C^{S}:S_{\operatorname{c}}\times S_{\operatorname{f}}\rightarrow S$ is given
by $C^{S}\left(  u_{\operatorname{c}},u_{\operatorname{f}}\right)
:=u_{\operatorname{c}}+u_{\operatorname{f}}$. For an operator $K:S\rightarrow
S$ and $s,t\in\left\{  \operatorname{c},\operatorname{f}\right\}  $ we set
$K_{st}:=\Pi_{s}^{S}K\Pi_{t}^{S}$.

Then, it follows from \cite[(3.17), Lem. 3.6]{grote2021stabilized}%
\begin{align*}
\left(  A^{S,p,\nu}\right)  ^{-1} &  =\left(  A^{S}\right)  ^{-1}+\left(
\left(  A^{S,p,\nu}\right)  ^{-1}-\left(  A^{S}\right)  ^{-1}\right)  \\
&  =\left(  A^{S}\right)  ^{-1}+C^{S}\left[
\begin{array}
[c]{cc}%
0 & 0\\
0 & \left(  A_{\operatorname{f}\operatorname{f}}^{S,p,\nu}\right)
^{-1}-\left(  A_{\operatorname{f}\operatorname{f}}^{S}\right)  ^{-1}%
\end{array}
\right]  E^{S}\\
&  =\left(  A^{S}\right)  ^{-1}+\left(  \left(  A_{\operatorname{f}%
\operatorname{f}}^{S,p,\nu}\right)  ^{-1}-\left(  A_{\operatorname{f}%
\operatorname{f}}^{S}\right)  ^{-1}\right)  \Pi_{\operatorname{f}}^{S}%
\end{align*}
and%
\begin{align}
\left(  A^{S,p,\nu}\right)  ^{-1}Q_{p,\nu,|r|}^{\Delta t}\left(  A^{S}%
\Pi_{\operatorname*{f}}^{S}\right)  \Pi_{\operatorname*{f}}^{S}  & =\left(
A^{S}\right)  ^{-1}Q_{p,\nu,|r|}^{\Delta t}\left(  A^{S}\Pi_{\operatorname*{f}%
}^{S}\right)  \Pi_{\operatorname*{f}}^{S}\label{estinv}\\
& +\left(  \left(  A_{\operatorname{f}\operatorname{f}}^{S,p,\nu}\right)
^{-1}-\left(  A_{\operatorname{f}\operatorname{f}}^{S}\right)  ^{-1}\right)
Q_{p,\nu,|r|}^{\Delta t}\left(  A_{\operatorname{f}\operatorname{f}}%
^{S}\right)  \Pi_{\operatorname*{f}}^{S},\nonumber
\end{align}
where we used that $\Pi_{\operatorname{f}}^{S}$ is a projection. Clearly
$A_{\operatorname{f}\operatorname{f}}^{S}$ and $Q_{p,\nu,|r|}^{\Delta
t}\left(  A_{\operatorname{f}\operatorname{f}}^{S}\right)  $ are self-adjoint
with respect to $\left(  \cdot,\cdot\right)  _{\mathcal{T}}$ and positive
semidefinite. Furthermore, \cite[Lem. 3.6]{grote2021stabilized} implies%
\[
\left\Vert \left(  A_{\operatorname{f}\operatorname{f}}^{S,p,\nu}\right)
^{-1}-\left(  A_{\operatorname{f}\operatorname{f}}^{S}\right)  ^{-1}%
\right\Vert _{S\leftarrow S}\leq\frac{\nu+1}{2\nu}\Delta t^{2}%
\]
and by applying well-known spectral theory we get%
\begin{equation}
\left\Vert \left(  \left(  A_{\operatorname{f}\operatorname{f}}^{S,p,\nu
}\right)  ^{-1}-\left(  A_{\operatorname{f}\operatorname{f}}^{S}\right)
^{-1}\right)  Q_{p,\nu,|r|}^{\Delta t}\left(  A_{\operatorname{f}%
\operatorname{f}}^{S}\right)  \Pi_{\operatorname*{f}}^{S}\right\Vert
_{\mathcal{T\leftarrow T}}\leq\frac{\nu+1}{2\nu}\Delta t^{2}\max
_{x\in\mathcal{\sigma}\left(  A_{\operatorname{f}\operatorname{f}}^{S}\right)
}\left\vert Q_{p,\nu,|r|}^{\Delta t}\left(  x\right)  \right\vert
.\label{estinvinv}%
\end{equation}
Next, we estimate the maximum on the right-hand side. The combination of
(\ref{eq:defQ_pnurk}) and (\ref{Def_beta}) with $Q_{p,\nu,r}^{\Delta
t}:=Q_{p,\nu,r,p}^{\Delta t}$ leads to%
\begin{equation}
Q_{p,\nu,\left\vert r\right\vert }^{\Delta t}(x)=\frac{p-\left\vert
r\right\vert }{p^{2}}\frac{\,U_{p-1-\left\vert r\right\vert }\left(
\delta_{p,\nu}-\frac{\Delta t^{2}x}{\omega_{p,\nu}}\right)  }%
{U_{p-1-\left\vert r\right\vert }\left(  \delta_{p,\nu}\right)  }%
.\label{Qexpl}%
\end{equation}
From \cite[(3.21)]{grote2021stabilized} it follows that the spectrum of
$A_{\operatorname{f}\operatorname{f}}^{S}$ is contained in $\left[
0,\frac{2\delta_{p,\nu}\omega_{p,\nu}}{\Delta t^{2}}\right]  $ so that the
argument $\delta_{p,\nu}-\frac{\Delta t^{2}x}{\omega_{p,\nu}}$ in
(\ref{Qexpl}) is contained in $\left[  -\delta_{p,\nu},\delta_{p,\nu}\right]
$. From \cite[18.7.3, 18.14.1]{NIST:DLMF} we conclude that%
\[
\left\Vert U_{p-1-\left\vert r\right\vert }\right\Vert _{L^{\infty}\left(
\left[  -1,1\right]  \right)  }=U_{p-1-\left\vert r\right\vert }\left(
1\right)
\]
and since $U_{p-1-\left\vert r\right\vert }$ behaves monotonically outside
$\left[  -1,1\right]  $ we obtain by symmetry%
\begin{equation}
\,\left\vert U_{p-1-\left\vert r\right\vert }\left(  \delta_{p,\nu}%
-\frac{\Delta t^{2}x}{\omega_{p,\nu}}\right)  \right\vert \leq
U_{p-1-\left\vert r\right\vert }\left(  \delta_{p,\nu}\right)  \quad\forall
x\in\mathcal{\sigma}\left(  A_{\operatorname{f}\operatorname{f}}^{S}\right)
.\label{estUpm1}%
\end{equation}
This leads to%
\[
\left\vert Q_{p,\nu,\left\vert r\right\vert }^{\Delta t}\left(  x\right)
\right\vert \leq\frac{p-\left\vert r\right\vert }{p^{2}}\quad\forall
x\in\mathcal{\sigma}\left(  A_{\operatorname{f}\operatorname{f}}^{S}\right)  .
\]
By inserting this into (\ref{estinvinv}) we get%
\[
\left\Vert \left(  \left(  A_{\operatorname{f}\operatorname{f}}^{S,p,\nu
}\right)  ^{-1}-\left(  A_{\operatorname{f}\operatorname{f}}^{S}\right)
^{-1}\right)  Q_{p,\nu,|r|}^{\Delta t}\left(  A_{\operatorname{f}%
\operatorname{f}}^{S}\right)  \Pi_{\operatorname*{f}}^{S}\right\Vert
_{\mathcal{T\leftarrow T}}\leq\frac{\nu+1}{2\nu}\Delta t^{2}\frac{p-\left\vert
r\right\vert }{p^{2}}.
\]
Next we estimate the first summand in (\ref{estinv}) and write%
\begin{equation}
\left(  A^{S}\right)  ^{-1}Q_{p,\nu,|r|}^{\Delta t}\left(  A^{S}%
\Pi_{\operatorname*{f}}^{S}\right)  \Pi_{\operatorname*{f}}^{S}=Q_{p,\nu
,|r|}^{\Delta t}\left(  0\right)  \left(  A^{S}\right)  ^{-1}\Pi
_{\operatorname*{f}}^{S}+\tilde{Q}_{p,\nu,|r|}^{\Delta t}\left(
A_{\operatorname{f}\operatorname{f}}^{S}\right)  \Pi_{\operatorname*{f}}%
^{S},\label{Asm1Q}%
\end{equation}
where%
\[
\tilde{Q}_{p,\nu,|r|}^{\Delta t}\left(  x\right)  :=\frac{Q_{p,\nu
,|r|}^{\Delta t}\left(  x\right)  -Q_{p,\nu,|r|}^{\Delta t}\left(  0\right)
}{x}.
\]
The expansion (\ref{firsttermexp2}) leads to%

\begin{align}
\left\Vert Q_{p,\nu,|r|}^{\Delta t}\left(  0\right)  \left(  A^{S}\right)
^{-1}\Pi_{\operatorname*{f}}^{S}\right\Vert _{S\leftarrow S}  & =\frac
{p-r}{p^{2}}\left\Vert \left(  A^{S}\right)  ^{-1}\Pi_{\operatorname*{f}}%
^{S}\right\Vert _{S\leftarrow S}\label{Asm1Q1stP}\\
& \overset{\text{\cite[(3.7)]{grote2021stabilized}}}{\leq}\frac
{C_{\operatorname*{eq}}^{2}}{c_{\operatorname{coer}}}\frac{p-r}{p^{2}%
}\left\Vert \Pi_{\operatorname*{f}}^{S}\right\Vert _{S\leftarrow S}\leq
\frac{C_{\operatorname*{eq}}^{2}}{c_{\operatorname{coer}}}\frac{p-r}{p^{2}%
}.\nonumber
\end{align}
For the estimate of $\tilde{Q}_{p,\nu,|r|}^{\Delta t}$ we employ (\ref{Qexpl})
for%
\[
\tilde{Q}_{p,\nu,|r|}^{\Delta t}\left(  x\right)  =\frac{p-\left\vert
r\right\vert }{p^{2}}\frac{\,U_{p-1-\left\vert r\right\vert }\left(
\delta_{p,\nu}-\frac{\Delta t^{2}x}{\omega_{p,\nu}}\right)  -U_{p-1-\left\vert
r\right\vert }\left(  \delta_{p,\nu}\right)  }{U_{p-1-\left\vert r\right\vert
}\left(  \delta_{p,\nu}\right)  x}.
\]
The combination of a Taylor argument and the monotonicity of orthogonal
polynomials outside $\left[  -1,1\right]  $ results in the estimate%
\[
\left\vert \tilde{Q}_{p,\nu,|r|}^{\Delta t}\left(  x\right)  \right\vert
\leq\frac{p-\left\vert r\right\vert }{p^{2}}\frac{U_{p-1-\left\vert
r\right\vert }^{\prime}\left(  \delta_{p,\nu}\right)  }{U_{p-1-\left\vert
r\right\vert }\left(  \delta_{p,\nu}\right)  }\frac{\Delta t^{2}}%
{\omega_{p,\nu}}\quad\forall x\in\mathcal{\sigma}\left(  A_{\operatorname{f}%
\operatorname{f}}^{S}\right)  .
\]
We employ Lemma \ref{LemA1} to estimate the ratio of the Chebyshev polynomials
and obtain%
\begin{equation}
\left\vert \tilde{Q}_{p,\nu,|r|}^{\Delta t}\left(  x\right)  \right\vert
\leq\frac{p-\left\vert r\right\vert }{p^{2}}\frac{\left(  p-\left\vert
r\right\vert \right)  ^{2}\operatorname*{e}^{\nu/2}}{3}\frac{\Delta t^{2}%
}{\omega_{p,\nu}}\quad\forall x\in\mathcal{\sigma}\left(  A_{\operatorname{f}%
\operatorname{f}}^{S}\right)  .\label{Asm1Q3rdP}%
\end{equation}
The estimate \cite[Lem. A.2]{grote2021stabilized} of $\omega_{p,\nu}$ implies%
\begin{equation}
\left\vert \tilde{Q}_{p,\nu,|r|}^{\Delta t}\left(  x\right)  \right\vert
\leq\frac{p-\left\vert r\right\vert }{p^{2}}\frac{\left(  p-\left\vert
r\right\vert \right)  ^{2}\operatorname*{e}^{\nu}}{6p^{2}}\Delta t^{2}%
\leq\frac{1}{6}\frac{p-\left\vert r\right\vert }{p^{2}}\operatorname*{e}%
\nolimits^{\nu}\Delta t^{2}.\label{Asm1Q2ndP}%
\end{equation}
This finishes the estimate of the left-hand side in (\ref{Asm1Q}) by combining
(\ref{Asm1Q1stP}) and (\ref{Asm1Q2ndP}):%
\[
\left\Vert \left(  A^{S}\right)  ^{-1}Q_{p,\nu,|r|}^{\Delta t}\left(  A^{S}%
\Pi_{\operatorname*{f}}^{S}\right)  \Pi_{\operatorname*{f}}^{S}\right\Vert
_{S\leftarrow S}\leq\frac{p-r}{p^{2}}\left(  \frac{C_{\operatorname*{eq}}^{2}%
}{c_{\operatorname{coer}}}+\frac{\operatorname*{e}\nolimits^{\nu}\Delta t^{2}%
}{6}\right)  .
\]
With this estimate at hand, the second term in (\ref{Rsl}) can be estimated
by
\begin{align}
&  \left\Vert \left(  A^{S,p,\nu}\right)  ^{-1}\sum_{r=-(p-1)}^{p-1}%
Q_{p,\nu,|r|}^{\Delta t}\left(  A^{S}\Pi_{\operatorname*{f}}^{S}\right)
\Pi_{\operatorname*{f}}^{S}\partial_{t}^{\ell}f_{S}\left(  t+\frac{r\Delta
t}{p}\right)  \right\Vert _{\mathcal{T}}\label{Qsum}\\
&  \qquad\leq\left(  \frac{C_{\operatorname*{eq}}^{2}}{c_{\operatorname{coer}%
}}+\frac{\operatorname*{e}\nolimits^{\nu}\Delta t^{2}}{6}\right)
\sum_{r=-(p-1)}^{p-1}\frac{p-\left\vert r\right\vert }{p^{2}}\max_{t\in\lbrack
t_{m-1},t_{m+2}]}\left\Vert \partial_{t}^{\ell}f_{S}\left(  t\right)
\right\Vert _{\mathcal{T}}\nonumber\\
&  \qquad=\left(  \frac{C_{\operatorname*{eq}}^{2}}{c_{\operatorname{coer}}%
}+\frac{\operatorname*{e}\nolimits^{\nu}\Delta t^{2}}{6}\right)  \max
_{t\in\lbrack t_{m-1},t_{m+2}]}\left\Vert \partial_{t}^{\ell}f_{S}\left(
t\right)  \right\Vert _{\mathcal{T}}.\nonumber
\end{align}

For the second estimate (\ref{lem:errorestRnEst2}) we use%
\begin{align}
\partial_{t}^{\ell}R_{S}\left(  t\right)   &  =P_{p,\nu}^{\Delta t}\left(
A^{S}\Pi_{\operatorname*{f}}^{S}\right)  \Pi_{\operatorname*{c}}^{S}%
f_{S}^{\left(  \ell\right)  }(t)+\left(  \sum_{r=-(p-1)}^{p-1}Q_{p,\nu
,|r|}^{\Delta t}\left(  A^{S}\Pi_{\operatorname*{f}}^{S}\right)  \right)
\Pi_{\operatorname*{f}}^{S}\partial_{t}^{\ell}f_{S}\left(  t\right)
\label{RsFdiff}\\
&  +\sum_{r=1}^{p-1}Q_{p,\nu,r}^{\Delta t}\left(  A^{S}\Pi_{\operatorname*{f}%
}^{S}\right)  \Pi_{\operatorname*{f}}^{S}\left(  \partial_{t}^{\ell}%
f_{S}\left(  t+\frac{r\Delta t}{p}\right)  +\partial_{t}^{\ell}f_{S}\left(
t-\frac{r\Delta t}{p}\right)  -2\partial_{t}^{\ell}f_{S}\left(  t\right)
\right)  .\nonumber
\end{align}
A Taylor argument for the last difference yields%
\[
\left\Vert \partial_{t}^{\ell}f_{S}\left(  t+\frac{r\Delta t}{p}\right)
+\partial_{t}^{\ell}f_{S}\left(  t-\frac{r\Delta t}{p}\right)  -2\partial
_{t}^{\ell}f_{S}\left(  t\right)  \right\Vert _{\mathcal{T}}\leq\sup
_{t\in\left[  t-\Delta t,t+\Delta t\right]  }\left\Vert \partial_{t}^{\ell
+2}f_{S}\left(  t\right)  \right\Vert _{\mathcal{T}}\left(  \frac{r\Delta
t}{p}\right)  ^{2}%
\]
and as in (\ref{Qsum}) we get%
\begin{align*}
& \left\Vert \left(  A^{S,p,\nu}\right)  ^{-1}\sum_{r=1}^{p-1}Q_{p,\nu
,|r|}^{\Delta t}\left(  A^{S}\Pi_{\operatorname*{f}}^{S}\right)
\Pi_{\operatorname*{f}}^{S}\left(  \partial_{t}^{\ell}f_{S}\left(
t+\frac{r\Delta t}{p}\right)  +\partial_{t}^{\ell}f_{S}\left(  t-\frac{r\Delta
t}{p}\right)  -2\partial_{t}^{\ell}f_{S}\left(  t\right)  \right)  \right\Vert
_{\mathcal{T}}\\
& \qquad\leq\Delta t^{2}\left(  \frac{C_{\operatorname*{eq}}^{2}%
}{c_{\operatorname{coer}}}+\frac{\operatorname*{e}\nolimits^{\nu}\Delta t^{2}%
}{6}\right)  \sup_{t\in\left[  t-\Delta t,t+\Delta t\right]  }\left\Vert
\partial_{t}^{\ell+2}f_{S}\left(  t\right)  \right\Vert _{\mathcal{T}}.
\end{align*}
For the first term in (\ref{RsFdiff}) we use (\ref{Ppnue2})%
\[
P_{p,\nu}^{\Delta t}\left(  A^{S}\Pi_{\operatorname*{f}}^{S}\right)
\Pi_{\operatorname*{c}}^{S}\partial_{t}^{\ell}f_{S}(t)=\Pi_{\operatorname*{c}%
}^{S}\partial_{t}^{\ell}f_{S}(t).
\]
For the second one we write%
\[
\sum_{r=-(p-1)}^{p-1}Q_{p,\nu,|r|}^{\Delta t}\left(  x\right)  =1+\rho\left(
x\right)
\]
for\
\[
\rho\left(  x\right)  :=\sum_{r=-(p-1)}^{p-1}\frac{p-\left\vert r\right\vert
}{p^{2}}\frac{\,U_{p-1-\left\vert r\right\vert }\left(  \delta_{p,\nu}%
-\frac{\Delta t^{2}x}{\omega_{p,\nu}}\right)  -U_{p-1-\left\vert r\right\vert
}\left(  \delta_{p,\nu}\right)  }{U_{p-1-\left\vert r\right\vert }\left(
\delta_{p,\nu}\right)  }.
\]
Then,%
\begin{align}
\left\Vert \left(  A^{S,p,\nu}\right)  ^{-1}\partial_{t}^{\ell}\left(
f_{S}-R_{S}\right)  \right\Vert _{\mathcal{T}}  & \leq\left\Vert \left(
A^{S,p,\nu}\right)  ^{-1}\partial_{t}^{\ell}\left(  f_{S}-\Pi
_{\operatorname*{c}}^{S}f_{S}-\Pi_{\operatorname*{f}}^{S}f_{S}\right)
\right\Vert _{\mathcal{T}}\label{mainsplit2ndterm}\\
& +\left\Vert \left(  A^{S,p,\nu}\right)  ^{-1}\rho\left(  A^{S}%
\Pi_{\operatorname*{f}}^{S}\right)  \Pi_{\operatorname*{f}}^{S}\partial
_{t}^{\ell}f_{S}\left(  t\right)  \right\Vert _{\mathcal{T}}\nonumber\\
& +\Delta t^{2}\left(  \frac{C_{\operatorname*{eq}}^{2}}%
{c_{\operatorname{coer}}}+\frac{\operatorname*{e}\nolimits^{\nu}\Delta t^{2}%
}{6}\right)  \sup_{t\in\left[  t-\Delta t,t+\Delta t\right]  }\left\Vert
\partial_{t}^{\ell+2}f_{S}\left(  t\right)  \right\Vert _{\mathcal{T}%
}.\nonumber
\end{align}
Clearly, the first term vanishes and it remains to estimate the second one.
Analogous arguments as for (\ref{estinv}) lead to%
\begin{equation}
\left(  A^{S,p,\nu}\right)  ^{-1}\rho\left(  A^{S}\Pi_{\operatorname*{f}}%
^{S}\right)  \Pi_{\operatorname*{f}}^{S}=\left(  A^{S}\right)  ^{-1}%
\rho\left(  A^{S}\Pi_{\operatorname*{f}}^{S}\right)  \Pi_{\operatorname*{f}%
}^{S}+\left(  \left(  A_{\operatorname{f}\operatorname{f}}^{S,p,\nu}\right)
^{-1}-\left(  A_{\operatorname{f}\operatorname{f}}^{S}\right)  ^{-1}\right)
\rho\left(  A_{\operatorname{f}\operatorname{f}}^{S}\right)  \Pi
_{\operatorname*{f}}^{S}\label{rhosplit}%
\end{equation}
and, in turn, as in (\ref{estinvinv}):%
\[
\left\Vert \left(  \left(  A_{\operatorname{f}\operatorname{f}}^{S,p,\nu
}\right)  ^{-1}-\left(  A_{\operatorname{f}\operatorname{f}}^{S}\right)
^{-1}\right)  \rho\left(  A_{\operatorname{f}\operatorname{f}}^{S}\right)
\Pi_{\operatorname*{f}}^{S}\right\Vert _{\mathcal{T\leftarrow T}}\leq\frac
{\nu+1}{2\nu}\Delta t^{2}\max_{x\in\mathcal{\sigma}\left(  A_{\operatorname{f}%
\operatorname{f}}^{S}\right)  }\left\vert \rho\left(  x\right)  \right\vert .
\]
It holds%
\[
\left\vert \rho\left(  x\right)  \right\vert \leq\sum_{r=-(p-1)}^{p-1}%
\frac{p-\left\vert r\right\vert }{p^{2}}\left\vert \frac{\,U_{p-1-\left\vert
r\right\vert }\left(  \delta_{p,\nu}-\frac{\Delta t^{2}x}{\omega_{p,\nu}%
}\right)  }{U_{p-1-\left\vert r\right\vert }\left(  \delta_{p,\nu}\right)
}-1\right\vert
\]
and as in (\ref{estUpm1}) we obtain%
\[
\left\vert \frac{\,U_{p-1-\left\vert r\right\vert }\left(  \delta_{p,\nu
}-\frac{\Delta t^{2}x}{\omega_{p,\nu}}\right)  }{U_{p-1-\left\vert
r\right\vert }\left(  \delta_{p,\nu}\right)  }\right\vert \leq1\qquad\forall
x\in\mathcal{\sigma}\left(  A_{\operatorname{f}\operatorname{f}}^{S}\right)  .
\]
The combination of these estimates implies%
\[
\left\Vert \left(  \left(  A_{\operatorname{f}\operatorname{f}}^{S,p,\nu
}\right)  ^{-1}-\left(  A_{\operatorname{f}\operatorname{f}}^{S}\right)
^{-1}\right)  \rho\left(  A_{\operatorname{f}\operatorname{f}}^{S}\right)
\Pi_{\operatorname*{f}}^{S}\right\Vert _{\mathcal{T\leftarrow T}}\leq\frac
{\nu+1}{\nu}\Delta t^{2}.
\]
It remains to estimate the first summand in the right-hand side of
(\ref{rhosplit}) and we proceed as in (\ref{Asm1Q}):%
\[
\left(  A^{S}\right)  ^{-1}\rho\left(  A^{S}\Pi_{\operatorname*{f}}%
^{S}\right)  \Pi_{\operatorname*{f}}^{S}=\rho\left(  0\right)  \left(
A^{S}\right)  ^{-1}\Pi_{\operatorname*{f}}^{S}+\tilde{\rho}\left(
A_{\operatorname{f}\operatorname{f}}^{S}\right)  \Pi_{\operatorname*{f}}^{S}%
\]
for%
\[
\tilde{\rho}\left(  x\right)  :=\frac{\rho\left(  x\right)  -\rho\left(
0\right)  }{x}.
\]
Since it holds by definition that $\rho\left(  0\right)  =0$, only the term containing
$\tilde{\rho}$ has to be estimated. It is easy to see that%
\[
\tilde{\rho}\left(  x\right)  =\frac{\rho\left(  x\right)  }{x}=\sum
_{r=-(p-1)}^{p-1}\tilde{Q}_{p,\nu,|r|}^{\Delta t}\left(  x\right)  .
\]
The estimate (\ref{Asm1Q2ndP}) then leads to%
\[
\left\Vert \tilde{\rho}\left(  A_{\operatorname{f}\operatorname{f}}%
^{S}\right)  \Pi_{\operatorname*{f}}^{S}\right\Vert _{\mathcal{T\leftarrow T}%
}\leq\sum_{r=-(p-1)}^{p-1}\frac{1}{6}\frac{p-\left\vert r\right\vert }{p^{2}%
}\operatorname*{e}\nolimits^{\nu}\Delta t^{2}=\frac{\operatorname*{e}%
\nolimits^{\nu}}{6}\Delta t^{2}.
\]
This finishes the proof of an estimate for the left-hand side in
(\ref{rhosplit}):%
\[
\left\Vert \left(  A^{S,p,\nu}\right)  ^{-1}\rho\left(  A^{S}\Pi
_{\operatorname*{f}}^{S}\right)  \Pi_{\operatorname*{f}}^{S}\right\Vert
_{\mathcal{T}\leftarrow\mathcal{T}}\leq\left(  \frac{\operatorname*{e}%
\nolimits^{\nu}}{6}+\frac{\nu+1}{\nu}\right)  \Delta t^{2}.
\]
By inserting this into (\ref{mainsplit2ndterm}) proves the second estimate
(\ref{lem:errorestRnEst2}):%
\[
\left\Vert \left(  A^{S,p,\nu}\right)  ^{-1}\partial_{t}^{\ell}\left(
f_{S}-R_{S}\right)  \right\Vert _{\mathcal{T}}\leq C_{3}\Delta t^{2}\max
_{k\in\left\{  0,2\right\}  }\sup_{t\in\left[  t-\Delta t,t+\Delta t\right]
}\left\Vert \partial_{t}^{\ell+k}f_{S}\left(  t\right)  \right\Vert
_{\mathcal{T}},
\]
for $C_{3}=\tfrac{C_{\operatorname*{eq}}^{2}}{c_{\operatorname{coer}}}%
+\tfrac{\nu+1}{\nu}+\tfrac{\operatorname*{e}\nolimits^{\nu}}{6}\left(
1+\Delta t^{2}\right)  $.%
\endproof

\begin{remark}
The ``split-LFC'' method proposed in \cite{CarleHochbruck3} for the inhomogeneous wave equation evaluates
the right-hand side only at each global time-step. 
The convergence analysis for such a coarser sampling of $f$ easily follows as a special case from of our theory: 
In Lemma \ref{lem:zSkn_explrepr}, the representation for that variant is obtained by setting $\Pi_{\operatorname{c}}^{S}$ to the identity operator and the last summation $\sum_{-k-1}^{k+1}\ldots$ to zero. 
Lemma \ref{lem:errorestRn} then holds verbatim while its proof already ends with estimate \eqref{eq:estPpnue2cnue}, since the analysis of the Chebyshev polynomials of second kind is only needed for a fully resolved right-hand side. 
As a consequence, the following two Theorems \ref{Theotimedisc} and \ref{TheoMain} also hold for that variant, that is, when $f$
is evaluated only once per global time-step, since the proofs remain identical.
\end{remark}

\begin{theorem}
\label{Theotimedisc} Assume that (\ref{wellposed}), (\ref{hsmallerh0}),
(\ref{ceqCeq}), (\ref{hfphc}) and (\ref{CFLtotal}) hold and let the solution
and right-hand side of the semi-discrete equation (\ref{spacedisc}) satisfy
$u_{S}\in W^{5,\infty}\left(  \left[  0,T\right]  ;L^{2}\left(  \Omega\right)
\right)  $ and $f_{S}\in W^{3,\infty}\left(  \left[  0,T\right]  ;L^{2}\left(
\Omega\right)  \right)  $, respectively. Then the fully discrete solution
$u_{S}^{\left(  n+1\right)  }$ of (\ref{eq:StabLFLTSrhs}) satisfies the error
estimate%
\[
\left\Vert e_{u,S,\Delta t}^{\left(  n+1\right)  }\right\Vert \leq C\Delta
t^{2}\left(  1+T\right)  \mathcal{M}\left(  u_{S},f_{S}\right)
\]
with
\[
\mathcal{M}\left(  u_{S},f_{S}\right)  :=\max\left\{  \max_{2\leq\ell\leq
5}\left\Vert \partial_{t}^{\ell}u_{S}\right\Vert _{L^{\infty}\left(  \left[
0,T\right]  ;L^{2}(\Omega)\right)  },\max_{0\leq\ell\leq3}\left\Vert
\partial_{t}^{\ell}f_{S}\right\Vert _{L^{\infty}\left(  \left[  0,T\right]
;L^{2}(\Omega)\right)  }\right\}
\]
and a constant $C$, which is independent of $n$, $\Delta t$, $T$, $h$, $p$,
$u_{S}$, and $f_{S}$.
\end{theorem}

\proof
The proof follows the lines of the proofs of \cite[Theorem 16]{grote_sauter_1}
and \cite[Theorem 3.10]{grote2021stabilized}. The LF-LTS$(\nu)$ scheme
\eqref{eq1} is stable \cite[Theorem 3.9]{grote2021stabilized}, i.e., for
$f_{S}\equiv0$,
\[
\left\Vert v_{S}^{\left(  n{+1/2}\right)  }\right\Vert _{\mathcal{T}%
}+\left\Vert u_{S}^{\left(  n\right)  }\right\Vert _{\mathcal{T}}\leq
C_{0}\left(  \left\Vert v_{S}^{\left(  1/2\right)  }\right\Vert _{\mathcal{T}%
}+\left\Vert u_{S}^{\left(  1\right)  }\right\Vert _{\mathcal{T}}\right)  ,
\]
where $C_{0}$ is independent of $n$, $\Delta t$, $h$, and $T$. We apply the
stability estimate to the second component of the error representation
(\ref{errorreprfu}) and obtain\footnote{For a pair of functions $\mathbf{v}%
=\left(  v_{1},v_{2}\right)  ^{\intercal}\in S\times S$ we use the notation
$\left\Vert \mathbf{v}\right\Vert _{\ell^{1}}:=\left\Vert v_{1}\right\Vert
_{\mathcal{T}}+\left\Vert v_{2}\right\Vert _{\mathcal{T}}$.}%
\begin{align}
\left\Vert e_{u,S,\Delta t}^{\left(  n+1\right)  }\right\Vert _{\mathcal{T}}
&  \leq C_{0}\left\Vert e_{S,\Delta t}^{\left(  1\right)  }\right\Vert
_{\mathcal{T}}+C_{0}\Delta t\sum_{\ell=1}^{n-1}\left\Vert
\mbox{\boldmath$ \sigma$}%
^{\left(  n-\ell\right)  }-%
\mbox{\boldmath$ \sigma$}%
^{\left(  n-\ell+1\right)  }\right\Vert _{\ell^{1}}\nonumber\\
&  +C_{0}\Delta t\left\Vert
\mbox{\boldmath$ \sigma$}%
^{\left(  1\right)  }\right\Vert _{\ell^{1}}+\Delta t\left\Vert
\mbox{\boldmath$ \sigma$}%
^{\left(  n\right)  }\right\Vert _{\ell^{1}}. \label{eestmainform}%
\end{align}

First, we consider the last term in the right-hand side of \eqref{eestmainform}.
By similar arguments as in the proof of \cite[Theorem 3.10]%
{grote2021stabilized}, we obtain
\begin{equation}
\Delta t\left\Vert \mbox{\boldmath$ \sigma$}^{\left(  n\right)  }\right\Vert
_{\ell^{1}}\leq\dfrac{\Delta t^{2}}{24}\left(  24\dfrac{\nu+1}{\nu}+\dfrac
{1}{c_{\nu}}+1\right)  \mathcal{M}_{n}\left(  u_{S},f_{S}\right)  +\left\Vert
\left(  A^{S,p,\nu}\right)  ^{-1}\left(  f_{S}^{(n)}-R_{S}^{(n)}\right)
\right\Vert _{\mathcal{T}}\label{estsigman}%
\end{equation}
with $c_{\nu}$ as in (\ref{Defcnue}) and%
\[
\mathcal{M}_{n}\left(  u_{S},f_{S}\right)  :=\max\left\{  \max_{2\leq\ell
\leq5}\left\{  \max_{t\in\left[  t_{n-1/2},t_{n+2}\right]  }\left\Vert
\partial_{t}^{\ell}u_{S}(t)\right\Vert _{\mathcal{T}}\right\}  ,\max
_{0\leq\ell\leq3}\left\{  \max_{t\in\left[  t_{n-1},t_{n+2}\right]
}\left\Vert \partial_{t}^{\ell}f_{S}(t)\right\Vert _{\mathcal{T}}\right\}
\right\}  .
\]
We apply Lemma \ref{lem:errorestRn} to \eqref{estsigman} to obtain
\begin{equation}
\Delta t\left\Vert \mbox{\boldmath$ \sigma$}^{\left(  n\right)  }\right\Vert
_{\ell^{1}}\leq\dfrac{\Delta t^{2}}{24}\left(  24\dfrac{\nu+1}{\nu}%
+\dfrac{1+C_{2}}{c_{\nu}}+1\right)  \mathcal{M}_{n}\left(  u_{S},f_{S}\right).
\end{equation}
Next, we investigate the summands in the second term of the right-hand side of
(\ref{eestmainform}),
\begin{align*}
&  \mbox{\boldmath$ \sigma$}^{\left(  m\right)  }%
-\mbox{\boldmath$ \sigma$}^{\left(  m+1\right)  }=\\
&  \!\!\left[  \!\!\!%
\begin{array}
[c]{c}%
\dfrac{u_{S}\left(  t_{m+2}\right)  -2u_{S}\left(  t_{m+1}\right)
+u_{S}\left(  t_{m}\right)  }{\Delta t^{2}}+\dfrac{v_{S}\left(  t_{m+1/2}%
\right)  -v_{S}\left(  t_{m+3/2}\right)  }{\Delta t}\\
\dfrac{u_{S}\left(  t_{m}\right)  -u_{S}\left(  t_{m+1}\right)  }{\Delta
t}-\left(  A^{S,p,\nu}\right)  ^{-1}\!\left(  \dfrac{v_{S}\left(
t_{m+3/2}\right)  -2v_{S}\left(  t_{m+1/2}\right)  +v_{S}\left(
t_{m-1/2}\right)  }{\Delta t^{2}}+\dfrac{R_{S}^{(m)}-R_{S}^{(m+1)}}{\Delta
t}\right)
\end{array}
\!\!\!\right]  ,
\end{align*}
where $R_S^{(m)}$ is defined as in Corollary \ref{CorLTS2step}.
By similar arguments as in the proof of \cite[Theorem 3.10]%
{grote2021stabilized} and by Lemma \ref{lem:errorestRn}, we obtain
\begin{align}
\left\Vert \mbox{\boldmath$ \sigma$}^{\left(  m\right)  }%
-\mbox{\boldmath$ \sigma$}^{\left(  m+1\right)  }\right\Vert _{\ell^{1}}%
&\leq\dfrac{\Delta t^{2}}{24}\left(  6\dfrac{\nu+1}{\nu}+\dfrac{3}{c_{\nu}%
}+4\right)  \mathcal{M}_{m}\left(  u_{S},f_{S}\right)  \nonumber\\
&\qquad+\left\Vert \left(
A^{S,p,\nu}\right)  ^{-1}\left(  \dot{f}_{S}^{(m+1/2)}-\dot{R}_{S}%
^{(m+1/2)}\right)  \right\Vert _{\mathcal{T}}.\label{estdiffsigmam}%
\end{align}
We apply Lemma\ref{lem:errorestRn} to \eqref{estdiffsigmam} to obtain
\begin{equation}
\left\Vert \mbox{\boldmath$ \sigma$}^{\left(  m\right)  }%
-\mbox{\boldmath$ \sigma$}^{\left(  m+1\right)  }\right\Vert _{\ell^{1}}%
\leq\dfrac{\Delta t^{2}}{24}\left(  6\dfrac{\nu+1}{\nu}+\dfrac{3+C_{2}}%
{c_{\nu}}+4\right)  \mathcal{M}_{m}\left(  u_{S},f_{S}\right)
\end{equation}
The rest of the proof follows by the same arguments as in the proof of
\cite[Theorem 3.10]{grote2021stabilized}.%
\endproof

\begin{theorem}
\label{TheoMain} Assume that (\ref{wellposed}), (\ref{hsmallerh0}),
(\ref{ceqCeq}), (\ref{hfphc}) and (\ref{CFLtotal}) hold and let the solution
and right-hand side of \eqref{waveeq} satisfy $u\in W^{8,\infty}\left(  \left[
0,T\right]  ;H^{m+1}\left(  \Omega\right)  \right)  $ and $f\in W^{6,\infty
}\left(  \left[  0,T\right]  ;H^{m+1}\left(  \Omega\right)  \right)  $,
respectively. Then, the corresponding fully discrete Galerkin FE formulation
with local time-stepping \eqref{eq:StabLFLTSrhs} has a unique solution which
satisfies the error estimate
\[
\left\Vert u\left(  t_{n+1}\right)  -u_{S}^{(n+1)}\right\Vert \leq
C(1+T)\left(  h^{m+1}+\Delta t^{2}\right)  \mathcal{Q}(u,f)
\]
with%
\[
Q(u,f):=\left(  1+Ch^{m+1}\left(  1+T\right)  \right)  \max\left\{  \left\Vert
u\right\Vert _{W^{8,\infty}\left(  \left[  0,T\right]  ;H^{m+1}\left(
\Omega\right)  \right)  },\left\Vert f\right\Vert _{W^{6,\infty}\left(
\left[  0,T\right]  ;H^{m+1}\left(  \Omega\right)  \right)  }\right\}
\]
and constants $C$ independent of $n$, $\Delta t$, $h$, $p$,
and the final time $T$.
\end{theorem}%

\proof
We adapt the proof in \cite[Thm .3.11]{grote2021stabilized} to the case of an inhomogeneous
right-hand side. We again use \cite[Thm 4.1]{Baker2}  to estimate the
error between the exact and semi-discrete solution. All assumptions in that
theorem are verified in the proof of \cite[Thm. 3.11]{grote2021stabilized} except for the
condition for the right-hand side. 
As usual, we denote below by $C$ a generic constant which may differ from one instance to the next.

From \cite[Lem. 5.2]{mass_lumping_2d}  it follows that%
\begin{equation}
\left\vert \left(  \chi,\psi\right)  -\left(  \chi,\psi\right)  _{\mathcal{T}%
}\right\vert \leq Ch^{r+\mu}\left(  \sum_{\tau\in\mathcal{T}}\left\Vert
\chi\right\Vert _{H^{r}\left(  \tau\right)  }^{2}\right)  ^{1/2}\left(
\sum_{\tau\in\mathcal{T}}\left\Vert \psi\right\Vert _{H^{\mu}\left(
\tau\right)  }^{2}\right)  ^{1/2}\quad\forall\chi,\psi\in S\label{chipiseest}%
\end{equation}
for all $1\leq r,\mu\leq m-1$. 
To deduce condition (3.5)
in \cite{Baker2} we transform the local quantities $\left\Vert
\chi\right\Vert _{H^{r}\left(  \tau\right)  }^{2}$ in (\ref{chipiseest}) to
the unit reference simplex $\hat{\tau}$. By $F_{\tau}:\hat{\tau}%
\rightarrow\tau$, we denote an affine bijection from the reference element to
the simplex $\tau$. For a function $g$ with domain $\tau,$ its pullback is
denoted by $\hat{g}_{\tau}=g\circ F_{\tau}$. Let $\left\vert \tau\right\vert $
denote the volume of the simplex $\tau\in\mathcal{T}$. A standard scaling
argument implies:%
\[
\left\Vert \chi\right\Vert _{H^{r}\left(  \tau\right)  }^{2}\leq C\left\vert
\tau\right\vert \sum_{\left\vert
\boldsymbol{\alpha}%
\right\vert \leq r}h_{\tau}^{-2\left\vert
\boldsymbol{\alpha}%
\right\vert }\left\Vert \partial^{%
\boldsymbol{\alpha}%
}\hat{\chi}\right\Vert _{L^{2}\left(  \hat{\tau}\right)  }^{2}.
\]
For $\mu\in\mathbb{N}$ and $q\geq2$, the embedding $W^{\mu,q}\left(  \hat
{\tau}\right)  \hookrightarrow L^{2}\left(  \hat{\tau}\right)  $ is continuous
so that%
\begin{align*}
\left\Vert \chi\right\Vert _{H^{r}\left(  \tau\right)  }^{2}  & \leq
C\left\vert \tau\right\vert \sum_{\left\vert
\boldsymbol{\alpha}%
\right\vert \leq r}h_{\tau}^{-2\left\vert
\boldsymbol{\alpha}%
\right\vert }\left\Vert \partial^{%
\boldsymbol{\alpha}%
}\hat{\chi}\right\Vert _{W^{\mu,q}\left(  \hat{\tau}\right)  }^{2}\\
& \leq C\left\vert \tau\right\vert \sum_{\left\vert
\boldsymbol{\alpha}%
\right\vert \leq r}h_{\tau}^{-2\left\vert
\boldsymbol{\alpha}%
\right\vert }\sum_{\left\vert
\boldsymbol{\beta}
\right\vert \leq\mu}\left\Vert \partial^{%
\boldsymbol{\alpha}%
+%
\boldsymbol{\beta}}\hat{\chi}\right\Vert _{L^{q}\left(  \hat{\tau}\right)  }^{2}.
\end{align*}
Transforming back to the physical simplex $\tau$ yields%
\[
\left\Vert \chi\right\Vert _{H^{r}\left(  \tau\right)  }^{2}\leq C\left\vert
\tau\right\vert \sum_{\left\vert
\boldsymbol{\alpha}%
\right\vert \leq r}h_{\tau}^{-2\left\vert
\boldsymbol{\alpha}%
\right\vert }\sum_{\left\vert
\boldsymbol{\beta}%
\right\vert \leq\mu}h^{2\left\vert
\boldsymbol{\alpha}%
+%
\boldsymbol{\beta}%
\right\vert }\left\vert \tau\right\vert ^{-1/q}\left\Vert \partial^{%
\boldsymbol{\alpha}%
+%
\boldsymbol{\beta}%
}\chi\right\Vert _{L^{q}\left(  \tau\right)  }^{2}.
\]
We use%
\[
h_{\tau}^{-2\left\vert
\boldsymbol{\alpha}%
\right\vert }h_{\tau}^{2\left\vert
\boldsymbol{\alpha}%
+%
\boldsymbol{\beta}%
\right\vert }\leq1\quad\text{and\quad}\left\vert \tau\right\vert \left\vert
\tau\right\vert ^{-2/q}\leq1
\]
to obtain%
\begin{equation}
\left\Vert \chi\right\Vert _{H^{r}\left(  \tau\right)  }^{2}\leq
C\sum_{\left\vert
\boldsymbol{\alpha}%
\right\vert \leq r}\sum_{\left\vert
\boldsymbol{\beta}%
\right\vert \leq\mu}\left\Vert \partial^{%
\boldsymbol{\alpha}%
+%
\boldsymbol{\beta}%
}\chi\right\Vert _{L^{q}\left(  \tau\right)  }^{2}\leq C\left\Vert
\chi\right\Vert _{W^{r+\mu,q}\left(  \tau\right)  }^{2}.\label{lastrefest}%
\end{equation}
Since the $\left\Vert \cdot\right\Vert _{\ell^{2}}-$ and the $\left\Vert
\cdot\right\Vert _{\ell^{q}}-$norm are equivalent in finite-dimensional
Euclidean (index) spaces, the combination of (\ref{chipiseest}) with
(\ref{lastrefest}) results in%
\begin{equation}
\left\vert \left(  \chi,\psi\right)  -\left(  \chi,\psi\right)  _{\mathcal{T}%
}\right\vert \leq Ch^{r+\mu}\left(  \sum_{\tau\in\mathcal{T}}\left\Vert
\chi\right\Vert _{W^{r+\mu,q}\left(  \tau\right)  }^{q}\right)  ^{1/q}\left(
\sum_{\tau\in\mathcal{T}}\left\Vert \psi\right\Vert _{H^{\mu}\left(
\tau\right)  }^{2}\right)  ^{1/2}\quad\forall\chi,\psi\in
S\label{lumpingerrorfem}%
\end{equation}
and this implies condition (3.5) in \cite[Thm 4.1]{Baker2}.

To estimate the effect of mass lumping in the right-hand side, let
$f:=f\left(  \cdot, t\right)  \in W^{r+\mu,q}\left(  \Omega\right)  $ for
$0<t\leq T$ and $\mu\in\left\{  1,2\right\}  $. We employ the
quasi-interpolation operator $\mathcal{I}_{h}^{\operatorname*{av}}%
:L^{1}\left(  \Omega\right)  \rightarrow S$ as in \cite[(5.1)]%
{Ern_Guermond_interpol} and the splitting:%
\begin{equation}
\left(  f,\psi\right)  -\left(  f,\psi\right)  _{\mathcal{T}}=\left(
f-\mathcal{I}_{h}^{\operatorname*{av}}f,\psi\right)  -\left(  f-\mathcal{I}%
_{h}^{\operatorname*{av}}f,\psi\right)  _{\mathcal{T}}+\left(  \mathcal{I}%
_{h}^{\operatorname*{av}}f,\psi\right)  -\left(  \mathcal{I}_{h}%
^{\operatorname*{av}}f,\psi\right)  _{\mathcal{T}}.\label{fsplit}%
\end{equation}
For the last difference we use (\ref{lumpingerrorfem}) to obtain%
\begin{equation}
\left\vert \left(  \mathcal{I}_{h}^{\operatorname*{av}}f,\psi\right)  -\left(
\mathcal{I}_{h}^{\operatorname*{av}}f,\psi\right)  _{\mathcal{T}}\right\vert
\leq Ch^{r+\mu}\left(  \sum_{\tau\in\mathcal{T}}\left\Vert \mathcal{I}%
_{h}^{\operatorname*{av}}f\right\Vert _{W^{r+\mu,q}\left(  \tau\right)  }%
^{q}\right)  ^{1/q}\left(  \sum_{\tau\in\mathcal{T}}\left\Vert \psi\right\Vert
_{H^{\mu}\left(  \tau\right)  }^{2}\right)  ^{1/2}.\label{Ihav1}%
\end{equation}
For the other terms in (\ref{fsplit}), the continuity of $\left(  \cdot
,\cdot\right)  $ and $\left(  \cdot,\cdot\right)  _{\mathcal{T}}$ yields%
\begin{equation}
\left\vert \left(  f-\mathcal{I}_{h}^{\operatorname*{av}}f,\psi\right)
\right\vert +\left\vert \left(  f-\mathcal{I}_{h}^{\operatorname*{av}}%
f,\psi\right)  _{\mathcal{T}}\right\vert \overset{\eqref{ceqCeq}}{\leq
}\left(  1+C_{\operatorname*{eq}}^{2}\right)  \left\Vert f-\mathcal{I}%
_{h}^{\operatorname*{av}}f\right\Vert \left\Vert \psi\right\Vert
.\label{Ihav2}%
\end{equation}
The approximation and stability properties of $\mathcal{I}_{h}%
^{\operatorname*{av}}$ are studied in \cite[Lem. 5.1, Thm. 5.2]%
{Ern_Guermond_interpol} and lead to%
\begin{equation}
\left\Vert f-\mathcal{I}_{h}^{\operatorname*{av}}f\right\Vert \leq Ch^{r+\mu
}\left\Vert f\right\Vert _{H^{r+\mu}\left(  \Omega\right)  }\label{Ihav3}%
\end{equation}
since $r+\mu\leq m+1$ and%
\begin{equation}
\left\Vert \mathcal{I}_{h}^{\operatorname*{av}}f\right\Vert _{W^{r+\mu
,q}\left(  \tau\right)  }\leq C\left\Vert f\right\Vert _{W^{r+\mu,q}\left(
\tau\right)  }.\label{Ihav4}%
\end{equation}
We combine (\ref{fsplit})--(\ref{Ihav4}) to deduce the estimate%
\[
\left\vert \left(  f,\psi\right)  -\left(  f,\psi\right)  _{\mathcal{T}%
}\right\vert \leq Ch^{r+\mu}\left(  \left\Vert f\right\Vert _{H^{r+\mu}\left(
\Omega\right)  }\left\Vert \psi\right\Vert +Ch^{r+\mu}\left\Vert f\right\Vert
_{W^{r+\mu,q}\left(  \Omega\right)  }\left(  \sum_{\tau\in\mathcal{T}%
}\left\Vert \psi\right\Vert _{H^{\mu}\left(  \tau\right)  }^{2}\right)
^{1/2}\right)  .
\]
Sobolev's embedding theorem implies $W^{r+\mu,q}\left(  \Omega\right)
\hookrightarrow H^{r+\mu}\left(  \Omega\right)  $ and $\left\Vert
\psi\right\Vert \leq\sqrt{\sum_{\tau\in\mathcal{T}}\left\Vert \psi\right\Vert
_{H^{\mu}\left(  \tau\right)  }^{2}}$ for $\mu\in\left\{  1,2\right\}  $ and,
finally, (3.8) in \cite[Thm 4.1]{Baker2} follows. The application of
this theorem yields the error estimate between the exact and semi-discrete
solution%
\begin{align}
\left\Vert \partial_{t}^{\ell}\left(  u-u_{S}\right)  \right\Vert _{L^{\infty
}\left(  \left[  0,T\right]  ;L^{2}\left(  \Omega\right)  \right)  }
&\leq
C_{\ell}h^{m+1}\left(  \sum_{k=\ell}^{\ell+1}\left\Vert \partial_{t}%
^{k}f\right\Vert _{L^{2}\left(  \left[  0,T\right]  ;H^{m+1}\left(
\Omega\right)  \right)  } \right.\nonumber\\
&\qquad\left. +\sum_{k=\ell}^{\ell+3}\left\Vert \partial_{t}%
^{k}u\right\Vert _{L^{2}\left(  \left[  0,T\right]  ;H^{m+1}\left(
\Omega\right)  \right)  }\right)  \label{estsemidisc}%
\end{align}
for $\ell=0$. Inspection of the proof in \cite[Thm 4.1]{Baker2}
implies that (\ref{estsemidisc}) also holds for any $\ell\in\mathbb{N}$
provided the right-hand side in (\ref{estsemidisc}) exists. Proceeding as in
the proof of \cite[Thm. 3.11]{grote2021stabilized} we arrive at%
\begin{align}
\left\Vert \partial_{t}^{\ell}\left(  u-u_{S}\right)  \right\Vert _{L^{\infty
}\left(  \left[  0,T\right]  ;L^{2}\left(  \Omega\right)  \right)  }
&\leq
C_{\ell}h^{m+1}\left(  1+T\right)  \left(  \left\Vert f\right\Vert
_{W^{\ell+1,\infty}\left(  \left[  0,T\right]  ;H^{m+1}\left(  \Omega\right)
\right)  } \right.\nonumber\\
&\qquad\left. +\left\Vert u\right\Vert _{W^{\ell+3,\infty}\left(  \left[
0,T\right]  ;H^{m+1}\left(  \Omega\right)  \right)  }\right)
.\label{deltumus}%
\end{align}
As in the proof of Theorem \cite[Thm. 3.11]{grote2021stabilized}, the combination with
(\ref{deltumus}) leads to%
\begin{align}
\left\Vert u\left(  t_{n+1}\right)  -u_{S}^{\left(  n+1\right)  }\right\Vert
&\leq C\left(  1+T\right)  \left(  h^{m+1}+\Delta t^{2}\right)  \max\left\{
\mathcal{M}\left(  u_{S},f_{S}\right)  , \right.\nonumber\\
&\quad\left. \left\Vert u\right\Vert _{W^{3,\infty
}\left(  \left[  0,T\right]  ;H^{m+1}\left(  \Omega\right)  \right)
},\left\Vert f\right\Vert _{W^{1,\infty}\left(  \left[  0,T\right]
;H^{m+1}\left(  \Omega\right)  \right)  }\right\}  .\label{semidiscr2nd}%
\end{align}
The final estimate is obtained by recalling the definition of
$\mathcal{M}$ as in Thm. \ref{Theotimedisc}:%
\begin{align*}
\left\Vert \partial_{t}^{\ell}u_{S}\right\Vert _{L^{\infty}\left(  \left[
0,T\right]  ;L^{2}\left(  \Omega\right)  \right)  }\leq & \left\Vert
\partial_{t}^{\ell}u\right\Vert _{L^{\infty}\left(  \left[  0,T\right]
;L^{2}\left(  \Omega\right)  \right)  }+\left\Vert \partial_{t}^{\ell}\left(
u_{S}-u\right)  \right\Vert _{L^{\infty}\left(  \left[  0,T\right]
;H^{m+1}\left(  \Omega\right)  \right)  }\\
\overset{\text{(\ref{deltumus})}}{\leq}  & \left\Vert \partial_{t}^{\ell
}u\right\Vert _{L^{\infty}\left(  \left[  0,T\right]  ;L^{2}\left(
\Omega\right)  \right)  }+\\
& +C_{\ell}h^{m+1}\left(  1+T\right)  \max\left\{  \left\Vert u\right\Vert
_{W^{\ell+3,\infty}\left(  \left[  0,T\right]  ;H^{m+1}\left(  \Omega\right)
\right)  }, \right. \\
&\quad\left. \left\Vert f\right\Vert _{W^{\ell+1,\infty}\left(  \left[
0,T\right]  ;H^{m+1}\left(  \Omega\right)  \right)  }\right\}  .
\end{align*}
Since $f\in W^{3,\infty}\left(  \left[  0,T\right]  ,L^{2}\left(
\Omega\right)  \right)  $, we obtain from the relation%
\[
\left(  \partial_{t}^{\ell}f_{S}\left(  t\right)  ,v\right)  =\left(
\partial_{t}^{\ell}f\left(  t\right)  ,v\right)  \quad\forall v\in S
\]
and upon setting $v=\partial_{t}^{\ell}f_{S}\left(  t\right)  \in S$:%
\[
\left\Vert \partial_{t}^{\ell}f_{S}\right\Vert _{L^{\infty}\left(  \left[
0,T\right]  ;L^{2}\left(  \Omega\right)  \right)  }\leq\left\Vert \partial
_{t}^{\ell}f\right\Vert _{L^{\infty}\left(  \left[  0,T\right]  ;L^{2}\left(
\Omega\right)  \right)  }.
\]
Hence,%
\[
\mathcal{M}\left(  u_{S},f_{S}\right)  \leq\left(  1+Ch^{m+1}\left(
1+T\right)  \right)  \max\left\{  \left\Vert u\right\Vert _{W^{8,\infty
}\left(  \left[  0,T\right]  ;H^{m+1}\left(  \Omega\right)  \right)
},\left\Vert f\right\Vert _{W^{6,\infty}\left(  \left[  0,T\right]
;H^{m+1}\left(  \Omega\right)  \right)  }\right\}
\]
and the combination with (\ref{semidiscr2nd}) leads to the assertion.%
\endproof

\section{Towards Optimal $H^{1}$ Convergence}\label{SecH1conv}
The LF-LTS($\nu$) Algorithm \ref{AlgStab} achieves the optimal convergence rate ${\cal O}(\Delta t^2 + h^{m+1})$
with respect to the $L^2$-norm with piecewise polynomial FE of degree $m$ also in the presence of a source $f$ -- see Theorem \ref{TheoMain}. In \cite{ChabassierImperiale2021}, however, it was shown that the original LF-LTS method from  \cite[Algorithm 1]{Grote_Mitkova}, that is, the LF-LTS($\nu$) Algorithm \ref{AlgStab} with $\nu=0$, may not
achieve the expected optimal convergence rate in the $H^1$-norm whenever the source $f$ is nonzero across the
coarse-to-fine interface. To restore optimal $H^1$-convergence regardless of $f$, we shall replace the somewhat abrupt projections from Section \ref{sec:FEM} by a new more gradual weighted transition. 

Let the set of simplex vertices be denoted by $\mathcal{V}\left(
\mathcal{T}\right)  $ and let $\mathcal{E}\left(
\mathcal{T}\right)  $ denote  the set of edges if $d\geq 2$ while for $d=1$ we set $\mathcal{E}\left(
\mathcal{T}\right)  = \mathcal{T}$. The definition of $S_{\mathcal{T}}^{m}$ gives rise to
the interpolation operator $I_{\mathcal{T}}^{m}:C^{0}\left(  \overline{\Omega
}\right)  \rightarrow S_{\mathcal{T}}^{m}$. For example, for $m=1$ the space
$S_{\mathcal{T}}^{1}$ consists of the continuous, piecewise linear polynomials
and $I_{\mathcal{T}}^{1}$ is the nodal interpolant%
\[
I_{\mathcal{T}}^{1}u=\sum_{z\in\mathcal{V}\left(  \mathcal{T}\right)
}u\left(  z\right)  b_{1,z}\text{.}%
\]
For two points $y,z$ in Euclidean space, we denote by $\left[  y,z\right]  $
the straight line connecting $y$ and $z$.

\begin{definition}
[connecting sequence]For $y,z\in\mathcal{V}\left(  \mathcal{T}\right)  $, a
sequence $\mathbf{x}=\left(  x_{i}\right)  _{i=0}^{m}\subset\mathcal{V}\left(
\mathcal{T}\right)  $ is a \emph{connecting sequence} of $y$ and $z$ if%
\[
x_{0}=y,\quad x_{m}=z,\quad\text{and for }1\leq j\leq m:\text{ }\left[
x_{j-1},x_{j}\right]  \in\mathcal{E}\left(  \mathcal{T}\right)  \text{.}%
\]
The length of $\mathbf{x}$ is $\operatorname{length}\left(  \mathbf{x}\right)
:=m$.
\end{definition}

\begin{definition}
[graph distance]The \emph{graph distance} between $y,z\in\mathcal{V}\left(
\mathcal{T}\right)  $ is given by%
\[
\operatorname{dist}_{\mathcal{T}}\left(  y,z\right)  :=\min\left\{
\operatorname{length}\left(  \mathbf{x}\right)  :\mathbf{x}\text{ is a
connecting sequence of }y\text{ and }z\right\}  .
\]
The \emph{graph distance} between $\mathcal{T}_{\operatorname{f}}$ and
some $z\in\mathcal{V}\left(  \mathcal{T}\right)  $ is given by%
\[
\operatorname{dist}\left(  \mathcal{T}_{\operatorname{f}},z\right)
:=\min_{y\in\mathcal{T}_{\operatorname{f}}}\operatorname{dist}_{\mathcal{T}%
}\left(  y,z\right)  .
\]

\end{definition}

\begin{definition}
[discrete distance function]Let $\mathcal{T}$ be a conforming simplicial finite element mesh
and let the subset $\mathcal{T}_{\operatorname{f}}\subset\mathcal{T}$ be
fixed. The \emph{discrete distance function for }$s$\emph{ simplex layers
about }$\mathcal{T}_{\operatorname{f}}$, $s\in\mathbb{N}_{\geq1}$, is given
as follows. For $z\in\mathcal{V}$, its value is set to%
\[
\eta_{s}\left(  z\right)  :=\left(  1-\frac{\operatorname{dist}\left(
\mathcal{T}_{\operatorname{f}},z\right)  }{s}\right)  _{+}%
\]
with $\left(  x\right)  _{+}:=\max\left\{  x,0\right\}  $. On $\Omega$, the
function is then given by
\[
\eta_{s}=\sum_{z\in\mathcal{V}\left(  \mathcal{T}\right)  }\eta_{s}\left(
z\right)  b_{1,z}.
\]

\end{definition}

\begin{definition}
[fine/coarse mesh mapping]The \emph{fine mesh mapping} $\Pi_{\operatorname{f}%
}^{S}:S\rightarrow S$ and \emph{coarse mesh mapping} $\Pi_{\operatorname{c}%
}^{S}:S\rightarrow S$ are given by%
\begin{equation}
\Pi_{\operatorname{f}}^{S}u:=\sum_{z\in\Sigma}u_{z}\eta_{s}\left(  z\right)
b_{z}\quad\text{and\quad}\Pi_{\operatorname{c}}^{S}:=I-\Pi_{\operatorname{f}%
}^{S}. \label{DefPismeared}%
\end{equation}

\end{definition}

\begin{remark}
\label{mapping}
Note that the above definition of $\Pi_{\operatorname{f}}^{S}$ and $\Pi
_{\operatorname{c}}^{S}$ generalizes the respective operators as defined in \cite{grote2021stabilized} and \eqref{defPifPic}: For $s=1$, both definitions coincide whereas for
$s>1$ the coarse-to-fine transition is now spread across several layers.
\end{remark}

Algorithm \ref{AlgStab} formally remains unchanged; however, the operators
$\Pi_{\operatorname{f}}^{S}$ and $\Pi_{\operatorname{c}}^{S}$ have to replaced
by the new definition (\ref{DefPismeared}). This comes at an extra cost since
the local time stepping is applied in $s$ additional simplex layers around
the original fine region $\Omega_{\operatorname*{f}}$ as in
(\ref{defsubdomains}) which are recursively defined by $\Omega
_{\operatorname*{f}}^{+1}$ as in (\ref{defsubdomains}) and for $r=2,\ldots s$
by%
\[
\Omega_{\operatorname*{f}}^{+r} 
:= \operatorname*{int} \left(  \bigcup {\left\{\tau\in\mathcal{T} : \tau\cap\overline{\Omega_{\operatorname*{f}}^{+\left(  r-1\right)  }}\neq\emptyset\right\}}\right)  .
\]

In Section \ref{SecH1convNumEx} we will report on numerical experiments for this modification.

\section{Numerical Experiments}
\label{SecNumEx}
First, we present numerical experiments which corroborate the convergence theory from Section \ref{SecConvAna}
and hence demonstrate that the LF-LTS($\nu$) method for the inhomogeneous wave equation defined in Algorithm \ref{AlgStab} achieves the expected optimal rates of convergence with respect to the $L^2$ norm. In doing so, 
we also compare the LF-LTS method with the alternative ``split-LFC'' approach from \cite{CarleHochbruck3}, where the
source inside the locally refined region is only evaluated once at each global time-step -- see Remark \ref{rem:LFCLTScomparison}.

Next, in Section 5.2, we show via two numerical experiments that the weighted transition, introduced in Section \ref{SecH1conv}, indeed restores optimal $H^1$-convergence regardless of the source and its position with respect to
the refined region in the mesh.

\subsection{Optimal $L^2$ Convergence}
\label{SecL2convNumEx}

We consider the one-dimensional wave equation \eqref{model problem} in $\Omega = (0,4)$ with 
homogeneous Dirichlet boundary conditions, i.e. $\Gamma = \Gamma_D$, $c \equiv 1$, and 
inhomogeneous source term 
\begin{equation}
f(x,t) = 250 \, e^{-400\left(\left(x-2\right)^2+\left(t-0.1\right)^2\right)},
\label{SourceEx1}
\end{equation}
a narrow space-time Gaussian pulse centered about $x=2$ and $t=0.1$.
The initial conditions are set to zero, i.e., $u_0 \equiv 0$, $v_0 \equiv 0$.

For the numerical solution, 
we use piecewise linear $H^1$-conforming finite elements with mass-lumping in space and the LF-LTS($\nu$) method 
with global time-step $\Delta t = \operatorname*{e}^{-\nu} h_{\operatorname*{c}}$ and $\nu = 0.01$ in time.
Next, we split the computational domain $\Omega$ into a coarse part, $\Omega_{\operatorname*{c}} = (0,1.6) \cup (2.4,4)$, and a locally refined part, $\Omega_{\operatorname*{f}} = [1.6,2.4]$, inside each of which we use an equidistant mesh with mesh sizes $h_{\operatorname*{c}}$ or $h_{\operatorname*{f}} = h_{\operatorname*{c}}/2$,  respectively.
Hence inside $\Omega_{\operatorname*{f}}$, the LF-LTS($\nu$) method takes $p=2$ local time-steps of size $\Delta t/p$ for each
global time-step of size $\Delta t$.

\begin{figure}
\centering
\includegraphics[width=\textwidth]{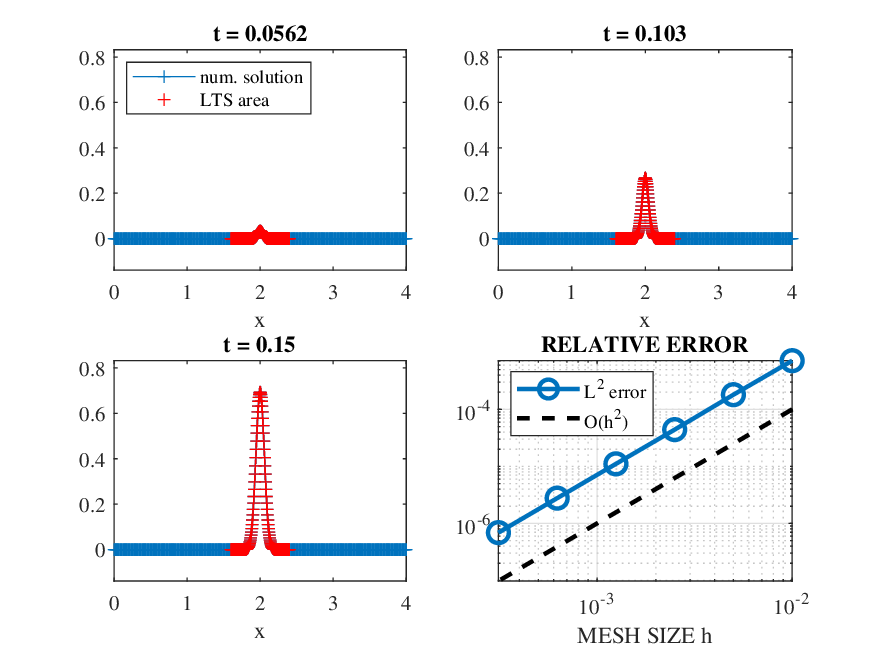}
\caption{Snapshots of the LF-LTS($\nu$) solution  of \eqref{model problem} with $\nu=0.01$ and $f$ as in \eqref{SourceEx1} (top left, top right and bottom left).
Relative $L^2$-error vs. $h = h_{\operatorname*{c}}$ for $\mathbb{P}_1$ finite elements (bottom right).}
\label{fig:solutions_L2conv}
\end{figure}

In Fig. \ref{fig:solutions_L2conv}, we display the LF-LTS($\nu$) numerical solution for $h_{\operatorname*{c}} = 0.01$ at times $t \approx 0.05$, $t \approx 0.1$, and $t = 0.15$, together with the relative $L^2$-error at $t=0.15$ for a sequence of locally refined meshes. The  $L^2$-error indeed displays 
the expected optimal second-order convergence rate $\mathcal{O}\left(h_{\operatorname*{c}}^2\right)$, as proved in Theorem \ref{TheoMain}.

\begin{figure}
\centering
\includegraphics[width=0.7\textwidth]{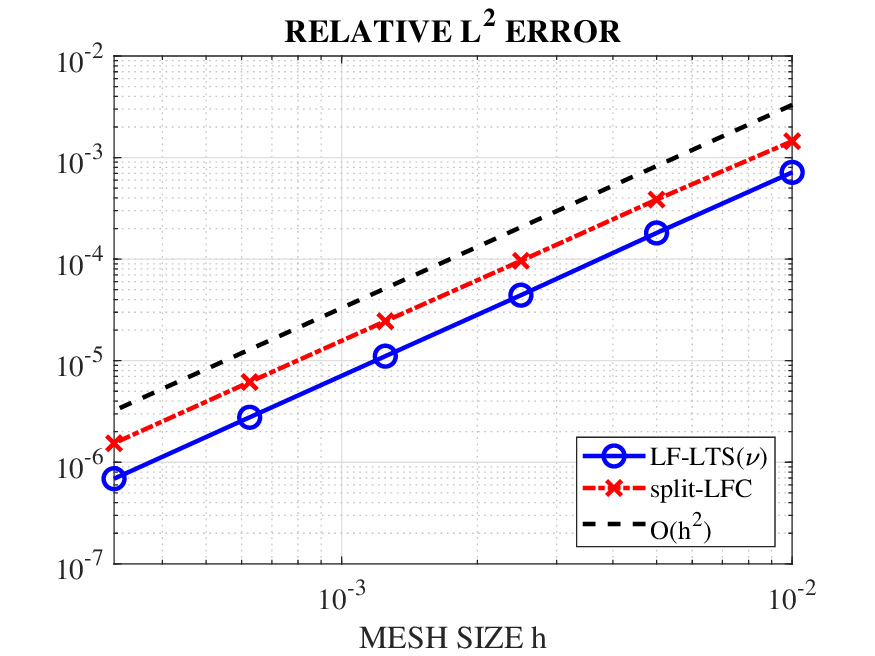}
\caption{$L^2$-error vs.  $h = h_{\operatorname*{c}}$: LF-LTS($\nu$) method \eqref{eq:StabLFLTSrhs} (dash-dotted red line) and  alternative ``split LFC'' approach from \cite{CarleHochbruck3} --  see Remark \ref{rem:LFCLTScomparison}  (solid blue line).}
\label{fig:compareL2errors}
\end{figure}

In Fig. \ref{fig:compareL2errors}, that same relative $L^2$-error is compared to that of the alternative approach
from \cite{CarleHochbruck3}, denoted here by LFC-LTS, which omits the intermediate source evaluations inside the refined region. 
Although both methods converge with the optimal rate $\mathcal{O}\left(h_{\operatorname*{c}}^2\right)$, as expected, the 
LF-LTS($\nu$)
algorithm from Section 2.3 here yields a threefold reduction in the error.

\subsection{Optimal $H^{1}$-Convergence and weighted transition}
\label{SecH1convNumEx}
In \cite{ChabassierImperiale2021}, it was shown for a particular (spatially constant) source term $f$ that the original LF-LTS method from  \cite[Algorithm 1]{Grote_Mitkova}, that is, the
 LF-LTS($\nu$) Algorithm \ref{AlgStab} with $\nu=0$, may not always achieve the expected optimal convergence rate in the $H^1$-norm.
Here, we carefully study this behavior for two different source terms and show that 
 the weighted transition from Section \ref{SecH1conv} restores the optimal $H^1$-convergence regardless of $f$.

\subsubsection{Space-time local forcing}

Again, we consider the wave equation \eqref{model problem} with zero initial conditions and source $f$ as in \eqref{SourceEx1} from the previous section.
For spatial discretization, however, we use piecewise quadratic $H^1$-conforming finite elements with mass-lumping ($\mathbb{P}_2$-FE).
The computational domain separates into a coarse part $\Omega_{\operatorname*{c}}$ and a locally refined part $\Omega_{\operatorname*{f}}$ with corresponding mesh sizes $h_{\operatorname*{c}}$ or $h_{\operatorname*{f}} = h_{\operatorname*{c}}/5$, respectively.

Now, while keeping the source fixed, we shall progressively shift the location of the refined region rightward and thus mimic the following three distinct
situations: (i) source $f$ located entirely inside $\Omega_{\operatorname*{f}}$, (ii) source $f$ nonzero across coarse-to-fine mesh interface, (iii)  source $f$ located entirely inside $\Omega_{\operatorname*{c}}$. In each case we monitor the $H^1$-convergence rates of the LF-LTS($\nu$) method for $p=5$ and $\nu=0.01$ either with, or without, weighted transition.
Hence the fine mesh mapping $\Pi^S_{\operatorname*{f}}$ is either defined by \eqref{DefPismeared}  with weighted transition across $s \sim 1/h_{\operatorname*{c}}$ elements, or without weighted transition, i.e., as in \eqref{DefPismeared} with $s=1$ -- see also Remark \ref{mapping}.

\begin{figure}
\centering
\includegraphics[width=0.45\textwidth]{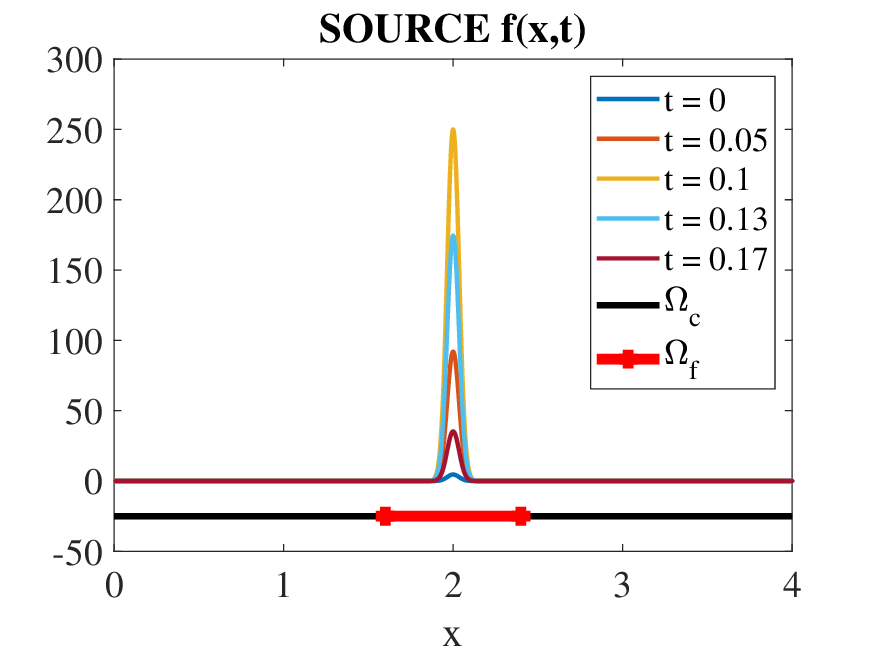}
\includegraphics[width=0.45\textwidth]{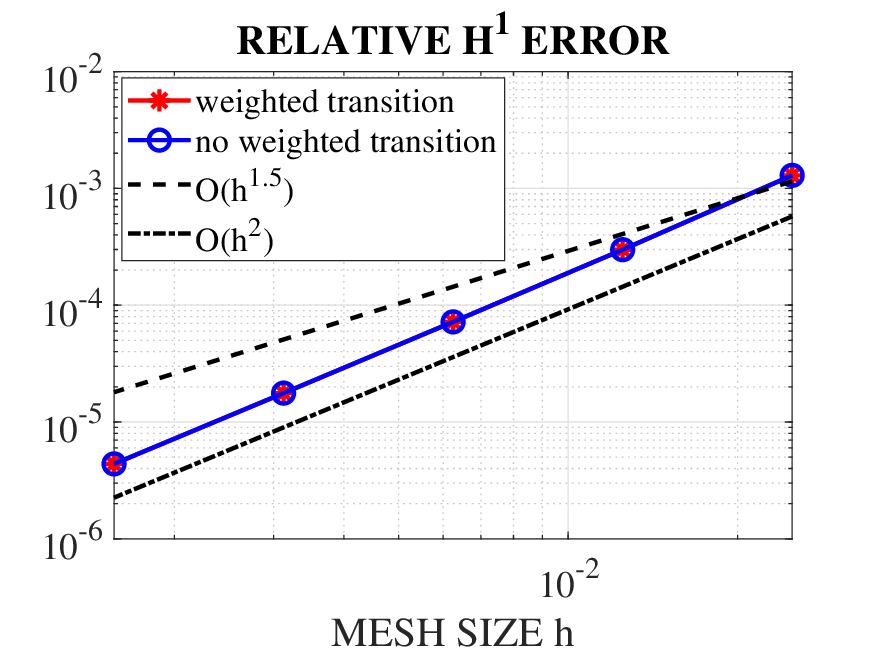}
\caption{Left: Source $f(x,t)$ for different times $t$, inside locally refined mesh with $\Omega_{\operatorname*{f}}=[1.6,2.4]$. Right: Relative $H^1$-error with (pink solid line), or without (blue dash-dotted line), weighted transition vs. mesh size $h = h_{\operatorname*{c}}$.}
\label{fig:H1conv_suppinside}
\end{figure}

First, we set $\Omega_{\operatorname*{f}} = [1.6,2.4]$ so that the support of $f$ completely lies inside $\Omega_{\operatorname*{f}}$ (up to machine precision), as illustrated
in Fig. \ref{fig:H1conv_suppinside}.
Here, the LF-LTS($\nu$) method achieves the optimal second-order rate $\mathcal{O}(h_{\operatorname*{c}}^2)$
 both with, or without, weighted transition, as shown in Fig. \ref{fig:H1conv_suppinside}. Recall that overall second-order
convergence vs. $h_{\operatorname*{c}}$ (or $\Delta t$) in the $H^1$-norm is expected here since we use
$\mathbb{P}_2$-FE in space. 

\begin{figure}
\centering
\includegraphics[width=0.45\textwidth]{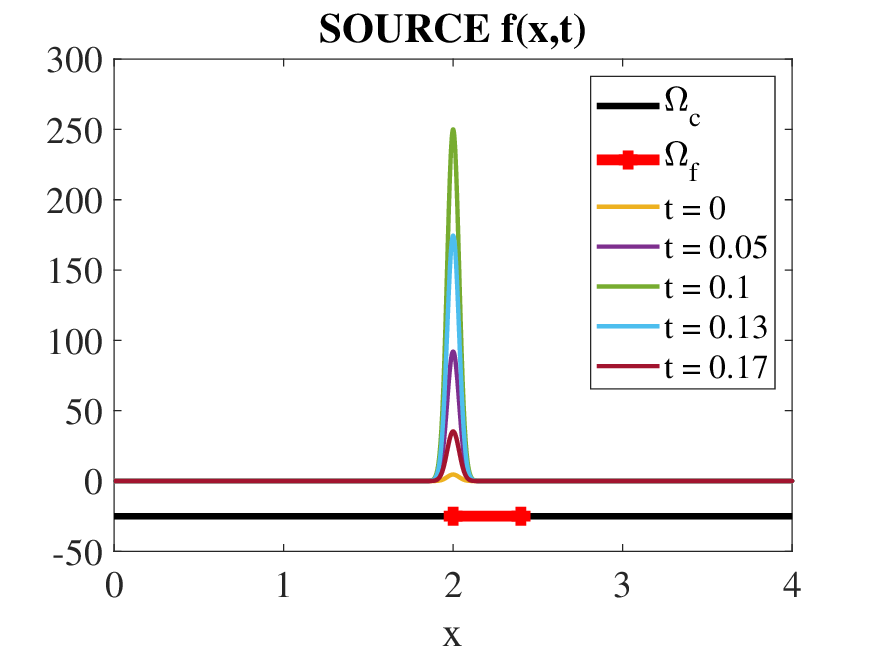}
\includegraphics[width=0.45\textwidth]{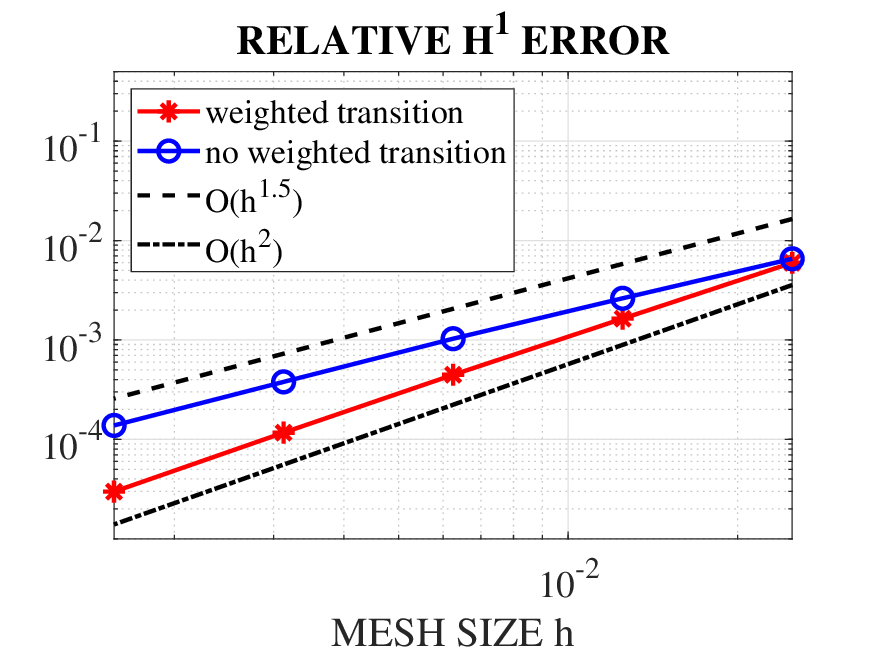} \caption{Left: Source $f(x,t)$ for different times $t$, nonzero across coarse-to-fine mesh interface with $\Omega_{\operatorname*{f}}=[2,2.4]$. Right: Relative $H^1$-error with (pink solid line), or without (blue dash-dotted line), weighted transition vs. mesh size $h = h_{\operatorname*{c}}$.}
\label{fig:H1conv_suppacross}
\end{figure}

Next, we shift the locally refined part rightwards and set $\Omega_{\operatorname*{f}} = [2,2.4]$ so that $f$ is non-zero across the interface between $\Omega_{\operatorname*{f}}$ and $\Omega_{\operatorname*{c}}$, as shown in Fig. \ref{fig:H1conv_suppacross};
in fact, the coarse-to-fine mesh interface now coincides with the maximum of $f$ at $x=2$.
In Fig. \ref{fig:H1conv_suppacross},
we observe that the LF-LTS($\nu$) method without weighted transition only converges as $\mathcal{O}(h_{\operatorname*{c}}^{1.5})$, while using the weighted transition across $s \sim 1/h_{\operatorname*{c}}$ layers of elements restores the expected optimal second-order convergence.

\begin{figure}
\centering
\includegraphics[width=0.45\textwidth]{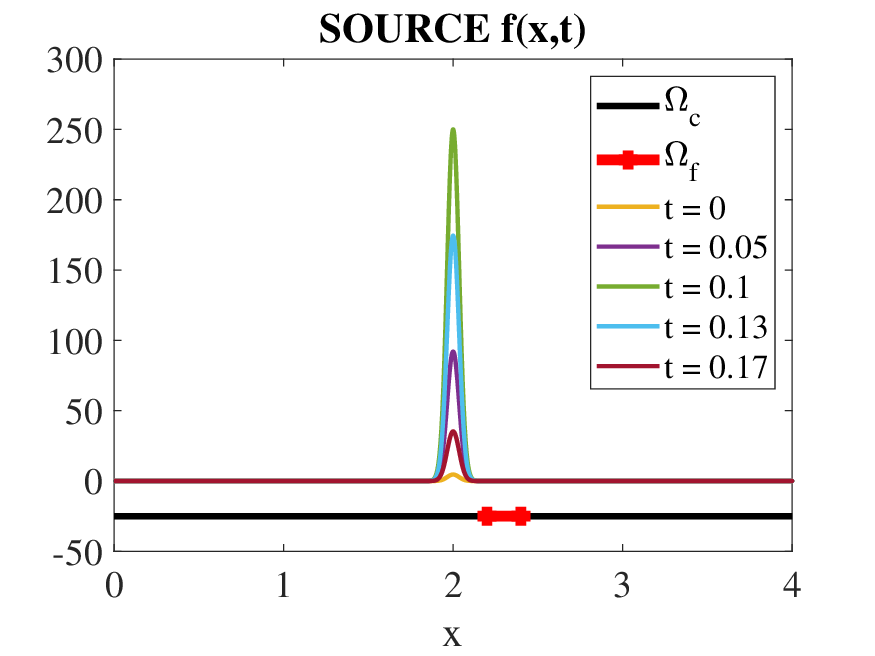}
\includegraphics[width=0.45\textwidth]{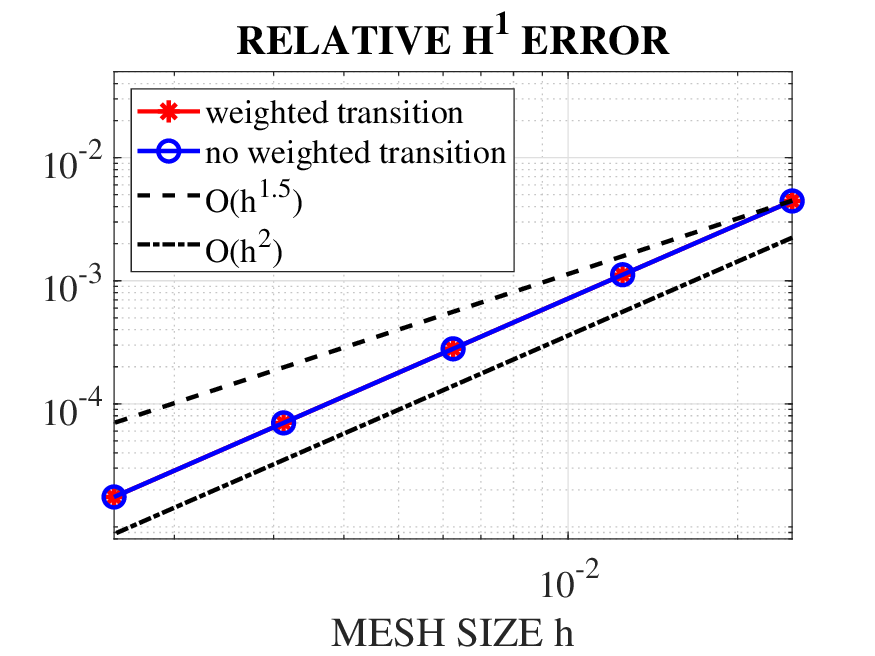}
\caption{Left: Source $f(x,t)$ for different times $t$, outside locally refined region with $\Omega_{\operatorname*{f}}=[2.2,2.4]$. Right: Relative $H^1$-error with (pink solid line), or without (blue dash-dotted line), weighted transition vs. mesh size $h = h_{\operatorname*{c}}$.}
\label{fig:H1conv_suppoutside}
\end{figure}

Finally, we shift  the locally refined part even farther to the right and set $\Omega_{\operatorname*{f}} := [2.2,2.4]$ so that the
 support of the source term essentially lies inside the coarse part of the mesh, as shown in Fig. \ref{fig:H1conv_suppoutside}.
Here, as in the first case, we again observe optimal second-order convergence in the $H^1$-norm, be it with or without weighted transition.

In summary, the LF-LTS($\nu$) method \emph{without weighted transition}
achieves in the considered examples the expected optimal convergence rates in the $H^1$-norm when the source vanishes
at the coarse-to-fine mesh transition; otherwise, when the source is non-zero across the mesh interface,
the $H^1$-convergence rate is reduced by one half. The \emph{weighted transition} introduced in Section 4 
restores the expected optimal $H^1$-convergence even when $f$ is nonzero across the coarse-to-fine
mesh interface.

\subsubsection{Spatially constant solution}

Finally, we revisit the example from \cite[Section 7.1.3]{ChabassierImperiale2021} and
consider the wave equation \eqref{model problem} in $\Omega = (0,1)$ with homogeneous Neumann boundary conditions
and zero initial conditions. Hence, we set the right-hand side $f$  as in  \cite{ChabassierImperiale2021} 
such that the solution is given by the spatially constant function
\begin{equation}
u(x,t) = \left\{
\begin{aligned}
&0, & &t \leq 0.1, \\
&\frac{1}{1 + \operatorname*{e}^{0.8\left(\frac{1}{t-0.1}+\frac{1}{t-0.9}\right)}}, & &0.1 < t < 0.9, \\
&1, & &t \geq 0.9.
\end{aligned}
\right.
\label{eq:ImperialeExe_sol}
\end{equation}

Now, we discretize \eqref{model problem} with $\mathbb{P}_2$-FE in space and the LF-LTS($\nu$) Algorithm \ref{AlgStab} with $\nu = 0.01$ in time. The refined region is (arbitrarily) set to the right half of the computational domain as $\Omega_{\operatorname*{f}} = [0.5,1]$ with a coarse-to-fine mesh size ratio of $p=2$. Clearly, here the spatially constant source is nonzero across 
the coarse-to-fine mesh transition.

\begin{figure}[t]
\centering
\includegraphics[width=0.7\linewidth]{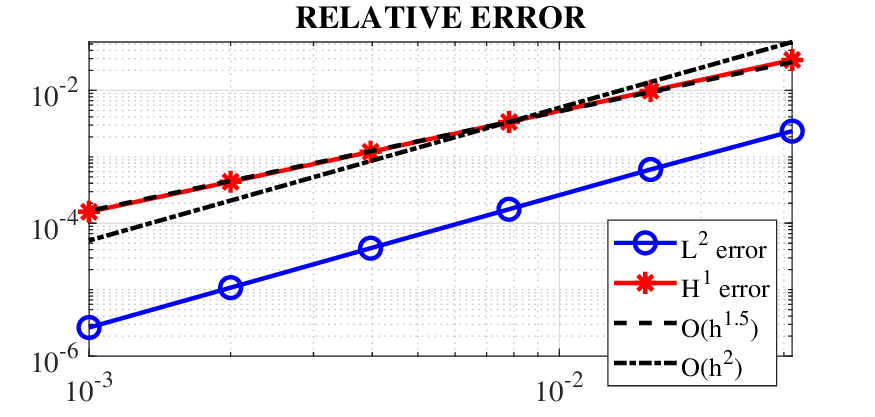}
\includegraphics[width=0.7\linewidth]{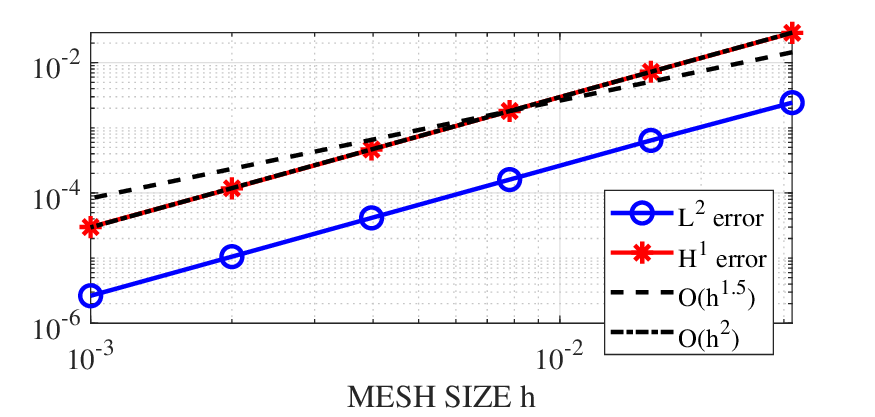}
\caption{Relative $L^2$- and $H^1$-errors without (top), or with (bottom), weighted transition.}
\label{fig:Convergence_weightedtrans}
\end{figure}
In Figure \ref{fig:Convergence_weightedtrans}, we compare the $L^2$- and $H^1$-errors on a sequence of meshes of size $h = h_{\operatorname*{c}}$ either with, or without, weighted transition.
Although the $L^2$-convergence theory from Section 3 does not apply here
due to the Neumann boundary conditions, we still observe optimal second-order convergence both with or without
weighted transition. Note that the overall second-order rate of convergence is dictated here by the time accuracy of the LF-LTS scheme. In contrast for the relative $H^1$-error, we observe an order reduction by $h^{1/2}$
without weighted transition. By using
the weighted transition across $s \sim 1/h$ element layers from Section 4, however,
optimal second-order convergence is again restored -- recall that $\mathbb{P}_2$-FE are used here
for the spatial discretization.

\begin{figure}[t]
\centering
\includegraphics[width=0.7\linewidth]{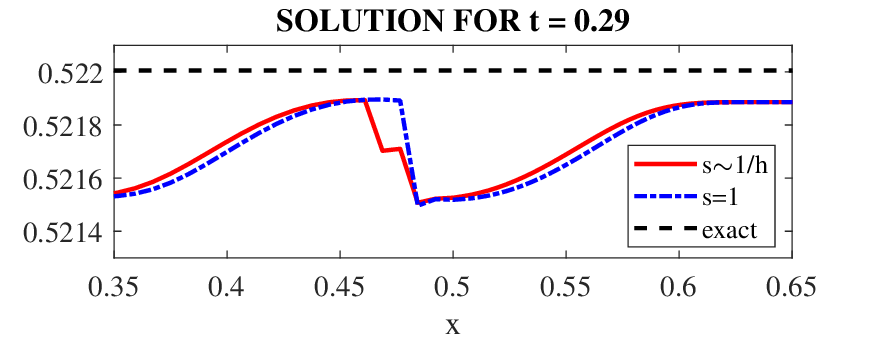}
\includegraphics[width=0.7\linewidth]{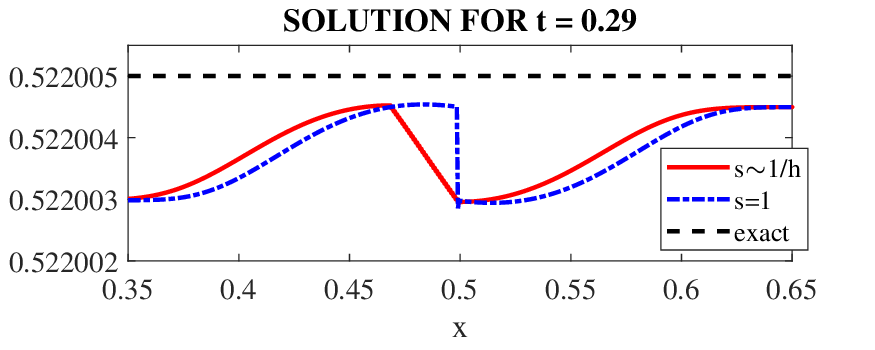}
\caption{Exact (dashed lines) and numerical solutions with ($s \sim 1/h$, solid), or without ($s=1$, dash-dotted), weighted transition at $t=0.29$ with $h=1.6\cdot10^{-2}$ (top) or $h=10^{-3}$ (bottom).}
\label{fig:Solutions_weightedtrans}
\end{figure}

To exhibit the mechanism that leads to this subtle order reduction by $h^{1/2}$,
which only affects the $H^1$ and not the $L^2$-convergence, we now strongly magnify in Fig. \ref{fig:Solutions_weightedtrans} the two (nearly constant) numerical solutions across the coarse-to-fine mesh transition for two different mesh sizes. 
Both variants of the LF-LTS($\nu$) method, with or without weighted transition, introduce a small perturbation near $x=0.5$, whose amplitude tends to zero as $h \rightarrow 0$. Without weighted transition, however, this perturbation increasingly steepens as $h \to 0$, which results in increased gradients at the intersection of $\Omega_{\operatorname*{c}}$ and $\Omega_{\operatorname*{f}}$.
Hence the weighted transition prevents the formation of increasingly steep gradients as $h\rightarrow 0$ 
at the coarse-to-fine mesh meeting point and thus restores the expected optimal $H^1$-convergence. 

\appendix

\section{Appendix}

Here we prove some technical estimates about Chebyshev polynomials.

\begin{lemma}
\label{LemA1}Let $\delta_{p,\nu}$ be as in (\ref{defdeltaomega}) and $0\leq
m\leq n\leq p-1$. Then%
\begin{equation}
\left\vert U_{n}^{\left(  m\right)  }\left(  \delta_{p,\nu}\right)
\right\vert \leq\frac{\left(  n+1\right)  ^{2m+1}\operatorname*{e}^{\nu/2}%
}{\left(  2m+1\right)  !!}.\label{Unm}%
\end{equation}
An estimate from below is given by%
\[
U_{n}^{\left(  m\right)  }\left(  \delta_{p,\nu}\right)  \geq2^{m}%
m!\binom{n+m+1}{n-m}\left(  \frac{2\nu}{p^{2}}\right)  ^{m}%
\]
with the particular case%
\[
U_{n}\left(  \delta_{p,\nu}\right)  \geq n+1.
\]

\end{lemma}%

\proof
From \cite[18.9.21]{NIST:DLMF}, we know that $U_{n}^{\left(  m\right)  }%
=\frac{1}{n+1}T_{n+1}^{\left(  m+1\right)  }$ and (\ref{Unm}) is equivalent to%
\[
T_{n+1}^{\left(  m+1\right)  }\left(  \delta_{p,\nu}\right)  \leq\frac{\left(
n+1\right)  ^{2m+2}\operatorname*{e}^{\nu/2}}{\left(  2m+1\right)  !!}.
\]
The definition $\delta_{p,\nu}=1+\nu/p^{2}$ along the side constraints for
$m,n,p-1$ imply%
\[
\delta_{p,\nu}\leq\delta_{n+1,\nu}.
\]
Since $T_{n+1}^{\left(  m+1\right)  }\left(  x\right)  $ is increasing for
$x\geq1$ we obtain from \cite[Lem. A.1]{grote2021stabilized}%
\[
T_{n+1}^{\left(  m+1\right)  }\left(  \delta_{p,\nu}\right)  \leq
T_{n+1}^{\left(  m+1\right)  }\left(  \delta_{n+1,\nu}\right)  \leq
\frac{\left(  n+1\right)  ^{2\left(  m+1\right)  }\operatorname*{e}^{\nu/2}%
}{\left(  2m+1\right)  !!}.
\]
For the estimate from below we use the representation in \cite[18.7.4 and
18.5.7]{NIST:DLMF}%

\begin{align*}
U_{n}^{\left(  m\right)  }\left(  x\right)    & =\frac{m!}{2^{m}}\frac
{n!}{\left(  \frac{3}{2}\right)  _{n}}\sum_{\ell=0}^{n}\frac{{\left(
n+2\right)  _{\ell}}{\left(  \ell+3/2\right)  _{n-\ell}}}{\ell!(n-\ell
)!}\binom{\ell}{m}\left(  \frac{x-1}{2}\right)  ^{\ell-m}\\
& =2^{m}m!\sum_{\ell=0}^{n}\binom{n+\ell+1}{n-\ell}\binom{\ell}{m}\left(
2\left(  x-1\right)  \right)  ^{\ell}.
\end{align*}
For $x=\delta_{p,\nu}$, this is a sum of positive terms and hence a lower
bound is obtained by truncating it after the first non-zero term, i.e., for
$\ell=m:$%
\[
U_{n}^{\left(  m\right)  }\left(  \delta_{p,\nu}\right)  \geq2^{m}%
m!\binom{n+m+1}{n-m}\left(  \frac{2\nu}{p^{2}}\right)  ^{m}.
\]%
\endproof

\bibliographystyle{abbrv}
\bibliography{references}

\begin{thebibliography}{10}

\bibitem{Baker2}
G.~A. Baker and V.~A. Dougalis.
\newblock The effect of quadrature errors on finite element approximations for
  second order hyperbolic equations.
\newblock {\em SIAM Journal on Numerical Analysis}, 13(4):577--598, 1976.

\bibitem{CarleHochbruck2}
C.~Carle and M.~Hochbruck.
\newblock Error analysis of multirate leapfrog-type methods for second-order
  semilinear {ODE}s.
\newblock {\em SIAM J. Numer. Anal.}, 60(5):2897--2924, 2022.

\bibitem{CarleHochbruck3}
C.~Carle and M.~Hochbruck.
\newblock Error analysis of second-order local time integration methods for
  discontinuous {G}alerkin discretizations of linear wave equations.
\newblock {\em Math. Comp.}, 93(350):2611--2641, 2024.

\bibitem{CarleHochbruchCheby}
C.~Carle, M.~Hochbruck, and A.~Sturm.
\newblock On leapfrog-{C}hebyshev schemes.
\newblock {\em SIAM J. Numer. Anal.}, 58(4):2404--2433, 2020.

\bibitem{ChabassierImperiale2015}
J.~Chabassier and S.~Imperiale.
\newblock Fourth-order energy-preserving locally implicit time discretization
  for linear wave equations.
\newblock {\em Internat. J. Numer. Methods Engrg.}, 106(8):593--622, 2016.

\bibitem{ChabassierImperiale2021}
{Chabassier, Juliette} and {Imperiale, Sébastien}.
\newblock Construction and convergence analysis of conservative second order
  local time discretisation for linear wave equations.
\newblock {\em ESAIM: M2AN}, 55(4):1507--1543, 2021.

\bibitem{mass_lumping_2d}
G.~Cohen, P.~Joly, J.~E. Roberts, and N.~Tordjman.
\newblock Higher order triangular finite elements with mass lumping for the
  wave equation.
\newblock {\em SIAM J. Numer. Anal.}, 38(6):2047--2078, 2001.

\bibitem{COHEN}
G.~C. Cohen.
\newblock {\em Higher {O}rder {N}umerical {M}ethods for {T}ransient {W}ave
  {E}quations}.
\newblock Springer Verlag, 2002.

\bibitem{ColFouJol3}
F.~Collino, T.~Fouquet, and P.~Joly.
\newblock A conservative space-time mesh refinement method for the 1-{D} wave
  equation. {I}. {C}onstruction.
\newblock {\em Numer. Math.}, 95(2):197--221, 2003.

\bibitem{ColFouJol2}
F.~Collino, T.~Fouquet, and P.~Joly.
\newblock A conservative space-time mesh refinement method for the 1-{D} wave
  equation. {II}. {A}nalysis.
\newblock {\em Numer. Math.}, 95(2):223--251, 2003.

\bibitem{ColFouJol1}
F.~Collino, T.~Fouquet, and P.~Joly.
\newblock Conservative space-time mesh refinement methods for the {FDTD}
  solution of {M}axwell's equations.
\newblock {\em J. Comput. Phys.}, 211(1):9--35, 2006.

\bibitem{DLM13}
S.~Descombes, S.~Lant\'eri, and L.~Moya.
\newblock Locally implicit discontinuous {G}alerkin method for time domain
  electromagnetics.
\newblock {\em J. Sci. Comp.}, 56:190--218, 2013.

\bibitem{DG09}
J.~Diaz and M.~J. Grote.
\newblock Energy conserving explicit local time stepping for second-order wave
  equations.
\newblock {\em SIAM J. Sci. Comput.}, 31(3):1985--2014, 2009.

\bibitem{DG15}
J.~Diaz and M.~J. Grote.
\newblock Multi-level explicit local time-stepping methods for second-order
  wave equations.
\newblock {\em Comput. Methods Appl. Mech. Engrg.}, 291:240--265, 2015.

\bibitem{NIST:DLMF}
{\it NIST Digital Library of Mathematical Functions}.
\newblock http://dlmf.nist.gov/, Release 1.0.13 of 2016-09-16.
\newblock F.~W.~J. Olver, A.~B. {Olde Daalhuis}, D.~W. Lozier, B.~I. Schneider,
  R.~F. Boisvert, C.~W. Clark, B.~R. Miller and B.~V. Saunders, eds.

\bibitem{DFFL10}
V.~Dolean, H.~Fahs, L.~Fezoui, and S.~Lanteri.
\newblock Locally implicit discontinuous {G}alerkin method for time domain
  electromagnetics.
\newblock {\em J.~Comput.~Phys.}, 229:512--526, 2010.

\bibitem{Ern_Guermond_interpol}
{Ern, Alexandre} and {Guermond, Jean-Luc}.
\newblock Finite element quasi-interpolation and best approximation.
\newblock {\em ESAIM: M2AN}, 51:1367--1385, 2017.

\bibitem{GilbertJoly}
J.~C. Gilbert and P.~Joly.
\newblock Higher order time stepping for second order hyperbolic problems and
  optimal {CFL} conditions.
\newblock In R.~Glowinski and P.~Neittaanm{\"a}ki, editors, {\em Partial
  Differential Equations}, volume~16 of {\em Comput. Methods Appl. Sci.}, pages
  67--93. Springer-Verlag, Dordrecht, The Netherlands, 2008.

\bibitem{grote_sauter_1}
M.~J. Grote, M.~Mehlin, and S.~A. Sauter.
\newblock Convergence analysis of energy conserving explicit local
  time-stepping methods for the wave equation.
\newblock {\em SIAM J. Numer. Anal.}, 56(2):994--1021, 2018.

\bibitem{grote2021stabilized}
M.~J. Grote, S.~Michel, and S.~A. Sauter.
\newblock Stabilized leapfrog based local time-stepping method for the wave
  equation.
\newblock {\em Math. Comp.}, 90(332):2603--2643, 2021.

\bibitem{Grote_Mitkova}
M.~J. Grote and T.~Mitkova.
\newblock Explicit local time-stepping methods for {M}axwell's equations.
\newblock {\em J. Comput. Appl. Math.}, 234(12):3283--3302, 2010.

\bibitem{GSS06}
M.~J. Grote, A.~Schneebeli, and D.~Sch{\"o}tzau.
\newblock Discontinuous {G}alerkin finite element method for the wave equation.
\newblock {\em SIAM J. Numer. Anal.}, 44(6):2408--2431, 2006.

\bibitem{HS16}
M.~Hochbruck and A.~Sturm.
\newblock Error analysis of a second-order locally implicit method for linear
  {M}axwell's equations.
\newblock {\em SIAM Journal on Numerical Analysis}, 54(5):3167--3191, 2016.

\bibitem{HS19}
M.~Hochbruck and A.~Sturm.
\newblock Upwind discontinuous {G}alerkin space discretization and locally
  implicit time integration for linear {M}axwell's equations.
\newblock {\em Math. Comp.}, 88(317):1121--1153, 2019.

\bibitem{HV03}
W.~Hundsdorfer and J.~Verwer.
\newblock {\em Numerical solution of time-dependent
  advection-diffusion-reaction equations}, volume~33 of {\em Springer Series in
  Computational Mathematics}.
\newblock Springer-Verlag, Berlin, 2003.

\bibitem{JolyRodriguez}
P.~Joly and J.~Rodr{\'{\i}}guez.
\newblock An error analysis of conservative space-time mesh refinement methods
  for the one-dimensional wave equation.
\newblock {\em SIAM J. Numer. Anal.}, 43(2):825--859, 2005.

\bibitem{LionsMagenesI}
J.~Lions and E.~Magenes.
\newblock {\em Non-{H}omogeneous {B}oundary {V}alue {P}roblems and
  {A}pplications}.
\newblock Springer-Verlag, Berlin, 1972.

\bibitem{MZKM13}
S.~Minisini, E.~Zhebel, A.~Kononov, and W.~A. Mulder.
\newblock Local time stepping with the discontinuous {G}alerkin method for wave
  propagation in 3{D} heterogeneous media.
\newblock {\em Geophysics}, 78:T67--T77, 2013.

\bibitem{MPFC08}
E.~Montseny, S.~Pernet, X.~Ferri\'{e}res, and G.~Cohen.
\newblock Dissipative terms and local time-stepping improvements in a spatial
  high order {D}iscontinuous {G}alerkin scheme for the time-domain {M}axwell's
  equations.
\newblock {\em J.~Comput.~Phys.}, 227:6795--6820, 2008.

\bibitem{Piperno}
S.~Piperno.
\newblock Symplectic local time-stepping in non-dissipative {DGTD} methods
  applied to wave propagation problems.
\newblock {\em M2AN Math. Model. Numer. Anal.}, 40(5):815--841, 2006.

\bibitem{Rietmann2017}
M.~Rietmann, M.~J. Grote, D.~Peter, and O.~Schenk.
\newblock Newmark local time stepping on high-performance computing
  architectures.
\newblock {\em J.~Comput.~Phys.}, 334:308--326, 2017.

\bibitem{Rivlin}
T.~J. Rivlin.
\newblock {\em The Chebyshev Polynomials}.
\newblock Wiley, New York, 1974.

\bibitem{Ver2011}
J.~G. Verwer.
\newblock Component splitting for semi-discrete {M}axwell equations.
\newblock {\em BIT}, 51:427--445, 2011.

\end{thebibliography}

\end{document}